\documentclass[preprint,12pt]{elsarticle}

\usepackage{amssymb}
\usepackage{amsmath}
\usepackage{bm}
\usepackage{amsthm}
\newtheorem{theorem}{Theorem}
\newtheorem{case}{Case}
\newtheorem{remark}{Remark}
\usepackage{tabularx}
\usepackage[ruled,vlined]{algorithm2e}
\usepackage{booktabs} 
\usepackage{subcaption}
\usepackage{float}

\begin{document}

\begin{frontmatter}



\title{ALM-PINNs Algorithms for Solving Nonlinear PDEs and Parameter Inversion Problems}




\author[1]{Yimeng Tian}
\ead{2023310133@stu.sufe.edu.cn}

\author[1]{Dinghua Xu\corref{cor1}}
\ead{dhxu6708@mail.shufe.edu.cn}

\address[1]{School of Mathematics, Shanghai University of Finance and Economics, Shanghai 200433, P.R.China}

\cortext[cor1]{Corresponding author}

\begin{abstract}
This paper focuses on the PINNs algorithm by proposing the ALM-PINNs computational framework to solve various nonlinear partial differential equations and corresponding parameters identification problems. The numerical solutions obtained by the ALM-PINNs algorithm are compared with both the exact solutions and the numerical solutions implemented from the PINNs algorithm. This demonstrates that under the same machine learning framework (TensorFlow 2.0) and neural network architecture, the ALM-PINNs algorithm achieves higher accuracy compared to the standard PINNs algorithm. Additionally, this paper systematically analyzes the construction principles of the loss function by introducing the probability distribution of random errors as prior information, and provides a theoretical basis for algorithm improvement.
\end{abstract}



\begin{keyword}
Physics-Informed Neural Networks(PINNs) \sep Nonlinear equations \sep Augmented Lagrangian Multiplier method(ALM) \sep Parameter Identification
\end{keyword}

\end{frontmatter}
\section{Introduction}
\label{sec:Intro}
In the past decade, rapid advancements have been made in artificial intelligence algorithms, with deep learning at the forefront, leading to widespread applications in key future industries such as integrated circuits, biopharmaceuticals, medical image analysis, and natural language processing. The computational research on neural networks has significantly accelerated the development of data-driven modeling approaches in the field of applied mathematics. While traditional numerical computing methods have been comprehensively and effectively applied to simple problems in low-dimensional regular domains, the development of numerical methods for high-dimensional nonlinear problems in more realistic and complex domains has been relatively slow.

\textbf{NNs for solving PDEs.} There is a gap between scientific research and industrial implementation. Weinan E.\cite{Weinan_2021} initiated and led research efforts in the development of deep learning-based algorithms in the fields of computational modeling and scientific computing. The development of neural networks, based on data-driven modeling, has provided new tools to tackle this challenge. Aristidis Likas\cite{Isaac_1998} first proposed a method using artificial neural networks to solve differential initial-boundary value problems. A. Pinkus\cite{Allan_1999} demonstrated that a sufficient number of neurons can spontaneously uniformly approximate any function and its partial derivatives, laying the groundwork for the development of solving differential equations with neural networks. Mohsen et al. \cite{Mohsen_2007} employed neural networks to solve the Burgers' equation within one-dimensional quasi-linear partial differential equations. Scientific computing tools such as automatic differentiation \cite{Baydin_2018} and its subsequent applications \cite{Kasa_2022} have further propelled the advancement of this field. 

Research into the neural networks solution of specific differential equation domains has been progressively deepening. For instance, Darbon et al. \cite{Darbon_2019} proposed a mathematical connection between Hamilton-Jacobi (HJ) partial differential equations with initial conditions and neural networks. Berner et al. \cite{Berner_2020} studied the numerical approximation of a specific class of d-dimensional Kolmogorov equations using neural networks. A. Jentzen et al. \cite{Jentzen_2021} demonstrated that deep artificial neural networks can overcome dimensionality issues in the numerical approximation of Kolmogorov partial differential equations with constant diffusion coefficients and nonlinear drift coefficients. 

\textbf{PINNs for solving PDEs and inverse problems.} The research direction of solving and estimating parameters of various differential equations with Physics-Informed Neural Networks (PINNs), has garnered wide attention among scholars. The ill-posedness of inverse problems often poses significant challenges to numerical computations. Different from traditional numerical integration algorithms for solving inverse problems, Raissi et al. \cite{Raissi_2019} first proposed the PINNs algorithm, which only requires the addition of extra terms in the construction of the objective functional of the forward problem. Its advantage lies in the comparability of the difficulty in solving both forward and inverse problems.

Hence, some research focuses on the direct application of PINNs to solve and estimate parameters of intricate differential equations. Raissi et al. \cite{Raissi_2020} introduced Hidden Fluid Mechanics (HFM): a deep learning framework based on Navier-Stokes information. Li et al. \cite{Li_2021} utilized PINNs to solve forward and inverse problems of the nonlinear Schrödinger equation. Shin et al. \cite{Shin_2020} investigated the convergence of Physics-Informed Neural Network solutions for linear second-order elliptic and parabolic equations. 
He et al. \cite{He_2020} discussed the use of PINNs algorithms for multi-layer physical data in underground transmission. Raissi et al. \cite{Raissi_2019} inverted parameters of Navier-Stokes equations in continuous-time models and Korteweg-de Vries equations in discrete-time models. PINNs can also be used to solve Integral Differential Equations (IDEs), Fractional Differential Equations (FDEs), and Stochastic Differential Equations (SDEs). Raissi et al. \cite{Raissi_2018}, Nabian et al. \cite{Na_2019}, Zhang et al. \cite{Zhang_2019}, Yang et al. have conducted research in this field. Lu et al. \cite{Lu_2021} developed the Python software package DeepXDE for computing various types of two-dimensional spatiotemporal problems. 

Some studies have improved the PINNs algorithm from a theoretical analysis perspective. Mojgani \cite{Mojgani_2023} introduced the concept of Lagrangian PINNs, Lu et al. \cite{Lu_2021} proposed a method based on DeepONet to address the issue of retraining the network for each new instance, Jagtap et al. \cite{Jagtap_2020} introduced adaptive activation functions into PINNs to accelerate the convergence speed of solutions. Bao et al. \cite{Bao_2020} investigated a class of weak adversarial neural networks for solving forward and inverse problems of partial differential equations. Basir et al. \cite{Basir_2022} studied artificial neural networks with equality constraints.

\textbf{Contribution of the paper:} This paper introduces the ALM-PINNs computational framework, an extension of the standard PINNs algorithm, to address nonlinear partial differential equations and parameter identification problems. The key contributions include demonstrating that ALM-PINNs achieve higher accuracy compared to PINNs under the same machine learning setup, and providing a systematic analysis of loss function construction, incorporating error distributions as prior information, which offers a theoretical basis for further algorithmic improvements.

\textbf{The structure of this paper.} Section \ref{sec:Intro} provides a brief introduction to the research directions and applications of artificial intelligence algorithms in solving forward and inverse differential equation problems. Section \ref{sec:ALM-PINNs-forward} presents the framework of the ALM-PINNs algorithm for solving forward problems. Section \ref{sec:ALM-PINNs-inverse} presents the framework of the ALM-PINNs algorithm for solving parameter inversion problems, derives the principles for constructing the error term in inverse problems, and proposes four strategies to address the ill-posedness in neural network training for inverse problems. Section \ref{sec:forward} provides two numerical examples for solving forward problems, comparing the numerical solutions of ALM-PINNs, PINNs, and exact solutions, which confirm that the ALM-PINNs algorithm achieves more accurate numerical results. Section \ref{sec:inverse} conducts parameter inversion for the numerical examples in the previous section, validating the theories proposed in Sections \ref{sec:ALM-PINNs-inverse}. Finally, Section \ref{sec:conclusion} presents the main conclusions of the paper.

\section{The construction principle of the ALM-PINNs algorithm}
\label{sec:ALM-PINNs-forward}
\subsection{Problem setup}
Let \( u(\bm{x}, t) \in H^2(0, T; H^2(\Omega)) \) be a function defined on the domain \( \Omega \times (0, t_T) \), where \( \Omega \) is an open region in \( \mathbb{R}^d \), and \( (0, t_T) \) is a time interval. The function \( u \) satisfies the following equation
\begin{equation}
    \label{gover_Eq}
    \mathcal{F}\left(\bm{x}, t; u, \frac{\partial u}{\partial t}, \frac{\partial^2 u}{\partial t^2}; \frac{\partial u}{\partial x_i}, \frac{\partial^2 u}{\partial x_i \partial x_j}, f,\bm{v}\right) = 0, \quad \bm{x} \in \Omega \subset \mathbb{R}^{d}, \quad t \in (0, t_T), \quad i,j=1,\ldots,d,
\end{equation}
subject to the initial boundary conditions respectively
\begin{equation}
    \mathcal{B}(\bm{x}, t; g, u) = u(\bm{x},t)-g(\bm{x},t) = 0, \quad \bm{x} \in \partial \Omega, \quad t \in [0, t_T], \quad i=1,\ldots,d,
\end{equation}
\begin{equation}
    \mathcal{I}(\bm{x}; u,h) = u(\bm{x},0)-h(\bm{x})=0, \quad \bm{x} \in \Omega, \quad i=1,\ldots,d,
\end{equation}
where $\mathcal{F}$ represents the residual form of the governing equation with differential operators, $\mathcal{B}$ represents the boundary conditions involving the known function $g$, and $\mathcal{I}$ represents the initial conditions involving the known function $h$. 
Specifically, $\mathcal{B}(\bm{x}, t; u, g)$ refers to the boundary condition, where $g$ is a given function that represents the prescribed boundary values, and $\mathcal{I}(\bm{x}; u, h)$ refers to the initial condition, where $h$ represents the known initial values of the solution.
Additionally, $f$ is the nonhomogeneous term of the equation, $\bm{x} = (x_1, \cdots, x_d)^T$ represents the $d$-dimensional spatial variables, $t$ is the time variable, $t_T$ is the terminal time, and $\bm{v}$ is the parameter vector.

\begin{figure}[htbp]
    \centering
    \includegraphics[width=0.5\linewidth]{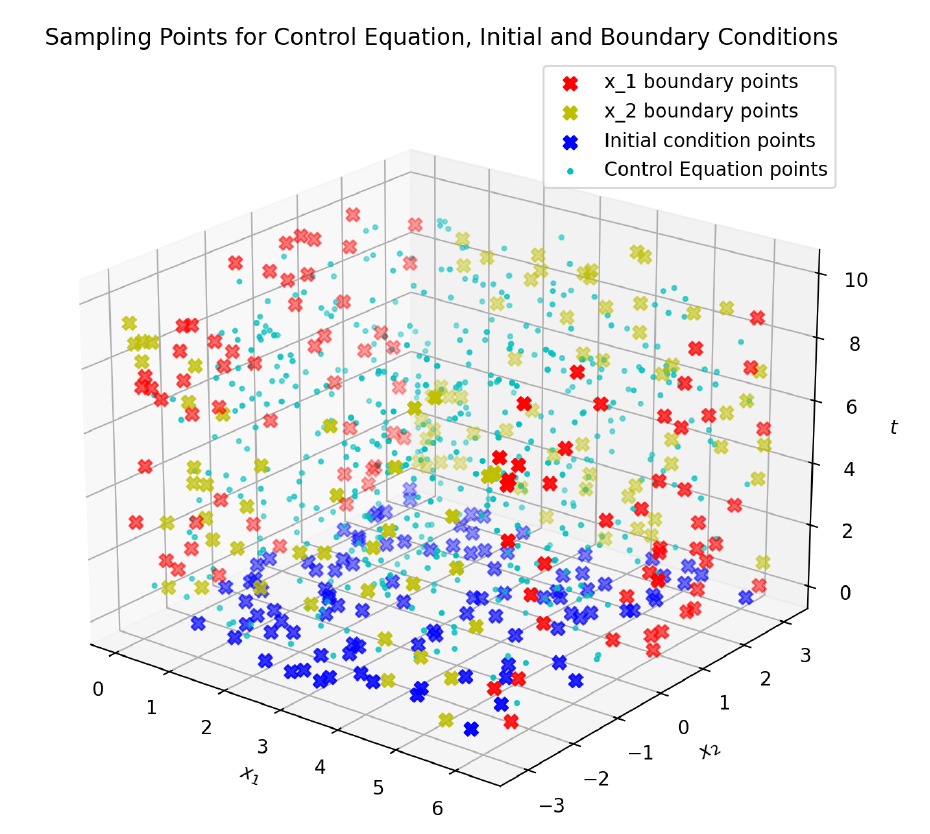}
    \caption{Sampling Points for Control Equation, Initial and Boundary Conditions}
    \label{fig:Samble}
\end{figure}

Training datasets $\{(\bm{x}^{(i)}_{inner}, t^{(i)}), f^{(i)}\}_{i=1}^{i=N_{\mathcal{F}}}$, $\{(\bm{x}^{(i)}_{boundary}, t^{(i)}), u^{(i)}, g^{(i)}\}_{i=1}^{i=N_{\mathcal{B}}}$, and $\{(\bm{x}^{(i)}_{inner}, 0), u^{(i)}, h^{(i)}\}_{i=1}^{i=N_{\mathcal{I}}}$ are given for the control equation and initial/boundary value conditions, respectively. Unlabeled data is provided for the control equation, while labeled data is available for initial and boundary conditions. $\lambda_{\mathcal{F}}$, $\lambda_{\mathcal{B}}$, and $\lambda_{\mathcal{I}}$ are hyper-parameters for the interior equation, boundary conditions, and initial conditions, respectively. The diagram in Figure \ref{fig:Samble} illustrates the sampling points for the internal, initial values, and boundary conditions of the equation when $d=2$.

Let \( u(\bm{x}, t; \bm{\theta}) \) represent the neural network solution, defined from the sample space \( \mathcal{X} \) to the output space \( \mathcal{Y} = \mathbb{R} \) or \( \mathbb{C} \). The solution is parameterized by \( \bm{\theta} \in \mathbb{R}^n \), a vector containing all the parameters of the neural network (i.e., weights and biases), where \( n \) denotes the total number of parameters.
The Hilbert space $H^2(0, T; H^2(\Omega))$ is assembled with the norm
\[
\|u\|_{H^2(0, T; H^2(\Omega))} = \left( \int_0^T \left( \|u(\cdot, t)\|_{H^2(\Omega)}^2 + \|\partial_t u(\cdot, t)\|_{H^2(\Omega)}^2 + \|\partial_{tt} u(\cdot, t)\|_{H^2(\Omega)}^2 \right) dt \right)^{1/2},
\]

The loss function for standard PINNs is defined as follows
\begin{equation}
    \mathcal{L}_{\mathcal{N}}(\bm{\theta}) = \lambda_{\mathcal{F}}\mathcal{L}_{\mathcal{F}}^{\mathcal{N}}(\bm{\theta})+\lambda_{\mathcal{B}}\mathcal{L}_{\mathcal{B}}^{\mathcal{N}}(\bm{\theta})+\lambda_{\mathcal{I}}\mathcal{L}_{\mathcal{I}}^{\mathcal{N}}(\bm{\theta}), 
    \label{eq:loss_PINNs}
\end{equation}

where the $\mathcal{L}_{\mathcal{F}}^{\mathcal{N}}(\bm{\theta})$ represents the loss function associated with the control equation component
\begin{equation}
    \mathcal{L}_{\mathcal{F}}^{\mathcal{N}}(\bm{\theta}) = \left\|\mathcal{F}(\bm{x}, t, f, u(\bm{x}, t;\bm{\theta})\right\|_{H^2(0, T; H^2(\Omega))}^2, 
    \label{eq:loss_control}
\end{equation}
the $\mathcal{L}_{\mathcal{B}}^{\mathcal{N}}(\bm{\theta})$ denotes the loss function corresponding to the boundary condition component
\begin{equation}
    \mathcal{L}_{\mathcal{B}}^{\mathcal{N}}(\bm{\theta}) = \left\|\mathcal{B}(\bm{x}, t, g, u(\bm{x}, t;\bm{\theta})\right\|_{H^2(0, T; H^2(\Omega))}^2,
    \label{eq:loss_boundary}
\end{equation}
and $\mathcal{L}_{\mathcal{I}}^{\mathcal{N}}(\bm{\theta})$ is the loss function for the initial value condition part
\begin{equation}
    \mathcal{L}_{\mathcal{I}}^{\mathcal{N}}(\bm{\theta}) = \left\|\mathcal{I}(\bm{x}, h, u(\bm{x};\bm{\theta})\right\|_{H^2(\Omega)}^2. 
\end{equation}

Note that function (\ref{eq:loss_PINNs}) can also be regarded as solving the following problem using the method of Lagrange multipliers. Considering the equation (\ref{eq:loss_PINNs}) as follows
\begin{equation}
    \mathcal{L}_{\mathcal{N}}(\bm{\theta}) = \mathcal{L}_{\mathcal{F}}^{\mathcal{N}}(\bm{\theta})+\frac{\lambda_{\mathcal{B}}}{\lambda_{\mathcal{F}}}\mathcal{L}_{\mathcal{B}}^{\mathcal{N}}(\bm{\theta})+\frac{\lambda_{\mathcal{I}}}{\lambda_{\mathcal{F}}}\mathcal{L}_{\mathcal{I}}^{\mathcal{N}}(\bm{\theta}), 
\end{equation}
consequently, the task of constructing the loss function can be formulated as the following optimization problem
\begin{equation}
\label{eq:primal}
    \begin{aligned} 
    &\min_{\bm{\theta}\in \mathbb{R}^n} \quad \mathcal{L}_{\mathcal{F}}^{\mathcal{N}}(\bm{\theta}),   \\
    &\begin{array}{r@{\quad}r@{}l@{\quad}l}
    s. t.   &    \mathcal{L}_{\mathcal{B}}^{\mathcal{N}}(\bm{\theta})=0, \\
    &   \mathcal{L}_{\mathcal{I}}^{\mathcal{N}}(\bm{\theta})=0. \\
    \end{array}
    \end{aligned}
\end{equation}

\subsection{Deterministic augmented Lagrangian method}
The Lagrangian function for the problem (\ref{eq:primal}) is
\begin{align}
    \mathcal{L}(\bm{\theta},\bm{\lambda}) &= \mathcal{L}_{\mathcal{F}}^{\mathcal{N}}(\bm{\theta}) + \lambda_B\mathcal{L}_{\mathcal{B}}^{\mathcal{N}}(\bm{\theta}) + \lambda_I\mathcal{L}_{\mathcal{I}}^{\mathcal{N}}(\bm{\theta}) \\
    &= \mathcal{L}_{\mathcal{F}}^{\mathcal{N}}(\bm{\theta}) + (\bm{\lambda}, \mathcal{L}_{\mathcal{C}}^{\mathcal{N}}(\bm{\theta})),
\end{align}
where $\bm{\lambda} = (\lambda_B, \lambda_I) \in \mathbb{R}^2$ is the Lagrange (dual) parameter associated with the constraint functions.
$\mathcal{L}_{\mathcal{C}}^{\mathcal{N}}(\bm{\theta}) = (\mathcal{L}_{\mathcal{B}}^{\mathcal{N}}(\bm{\theta}),\mathcal{L}_{\mathcal{I}}^{\mathcal{N}}(\bm{\theta}))^T$ is the matrix of initial-boundary conditions. Noting that
\begin{equation}
    \sup_{\bm{\lambda} \in \mathbb{R}^2} \mathcal{L}(\bm{\theta},\bm{\lambda}) =
    \left\{
    \begin{aligned}
    & \mathcal{L}_{\mathcal{F}}^{\mathcal{N}}(\bm{\theta}), &&  \mathcal{L}_{\mathcal{C}}^{\mathcal{N}}(\bm{\theta}) = 0, \\
    & \infty,    &&  \mathcal{L}_{\mathcal{C}}^{\mathcal{N}}(\bm{\theta}) \neq 0.
    \end{aligned}
    \right.
\end{equation}
Thus, if there exists $\bm{\theta}$ such that $\mathcal{L}_{\mathcal{C}}^{\mathcal{N}}(\bm{\theta}) = 0 $, 
\begin{equation}
    \min_{\bm{\theta} \in \mathbb{R}^n} \sup_{\bm{\lambda} \in \mathbb{R}^2} 
    \mathcal{L}(\bm{\theta},\bm{\lambda})
    = \min_{\{ \bm{\theta} \in \mathbb{R}^n |\mathcal{L}_{\mathcal{C}}^{\mathcal{N}}(\bm{\theta}) = 0 \}} \mathcal{L}_{\mathcal{F}}^{\mathcal{N}}(\bm{\theta}),
\end{equation}
and
\begin{equation}
    \arg \min_{\bm{\theta}\in \mathbb{R}^n} \sup_{\bm{\lambda} \in \mathbb{R}^2} 
    \mathcal{L}(\bm{\theta},\bm{\lambda})
    = \arg \min_{\{ \bm{\theta} \in \mathbb{R}^n |\mathcal{L}_{\mathcal{C}}^{\mathcal{N}}(\bm{\theta})= 0 \}} \mathcal{L}_{\mathcal{F}}^{\mathcal{N}}(\bm{\theta}).
\end{equation}
We define the saddle-point problem as
\begin{equation}
    \label{eq:dual}
    \min_{\bm{\theta} \in \mathbb{R}^n} \sup_{\bm{\lambda} \in \mathbb{R}^2} \mathcal{L}(\bm{\theta},\bm{\lambda}),
\end{equation}
which is equivalent to (\ref{eq:primal}).

A primal-dual iterate pair $(\bm{\theta}^*,\bm{\lambda}^*)$ is said to be a stationary point of (\ref{eq:dual}) if 
\begin{equation}
\label{eq:stationary_point}
    (\bm{\theta}^*,\bm{\lambda}^*) \in \{(\bm{\theta},\bm{\lambda})|\frac{P(\bm{\theta} - \eta {\nabla}_{\bm{\theta}}\mathcal{L}(\bm{\theta},\bm{\lambda}))-\bm{\theta}}{\eta} = 0 \quad and \quad \mathcal{L}_{\mathcal{C}}^{\mathcal{N}}(\bm{\theta}) = 0\},
\end{equation}
where $\eta > 0 $ and $P(\bm{x}) = \arg \min_{\bm{\theta}} \|\bm{\theta}-\bm{x}\|^2$ is the projection of $x \in \mathbb{R}^n$ onto the set $\bm{\theta}$. The condition that the gradient with respect to the projection equals zero indicates that at $\bm{\theta}$, the update along the gradient direction does not alter $\bm{\theta}$, ensuring it remains unchanged.

The augmented Lagrangian method is an iterative approach that generates stationary points satisfying (\ref{eq:stationary_point}) by solving a series of sub-problems. In each sub-problem, the objective function consists of the Lagrangian function $\mathcal{L}(\bm{\theta},\bm{\lambda})$ augmented by a quadratic penalty term, which penalizes violations of the equality constraint $\mathcal{L}_{\mathcal{C}}^{\mathcal{N}}(\bm{\theta}) = (\mathcal{L}_{\mathcal{B}}^{\mathcal{N}}(\bm{\theta}),\mathcal{L}_{\mathcal{I}}^{\mathcal{N}}(\bm{\theta}))^T = \bm{0}$. Specifically, at each iteration $k \in \mathbb{Z}^+$, the primal and dual update rules are defined as follows
\begin{equation}
    \bm{\theta}^*_k \in \arg \min_{\bm{\theta}} \mathcal{L}(\bm{\theta},\bm{\lambda};\mu_k),
\end{equation}
\begin{equation}
    \bm{\lambda}_{k+1} = \bm{\lambda}_k + \mu_k \mathcal{L}_{\mathcal{C}}^{\mathcal{N}}(\bm{\theta}_k^*),
\end{equation}
where $\mu_k >0 $ is the penalty parameter and
\begin{align}
\label{eq:ALM_Forward}
    \mathcal{L}_{Forward}(\bm{\theta};\bm{\lambda}, \bm{\mu}) 
    &= 
    \mathcal{L}_{\mathcal{F}}^{\mathcal{N}}(\bm{\theta})+
    (\bm{\lambda},\mathcal{L}_{\mathcal{C}}^{\mathcal{N}}(\bm{\theta}))+\frac{\mu_{\mathcal{B}}}{2}|\mathcal{L}_{\mathcal{B}}^{\mathcal{N}}(\bm{\theta})|^2+\frac{\mu_{\mathcal{B}}}{2}|\mathcal{L}_{\mathcal{I}}^{\mathcal{N}}(\bm{\theta})|^2, 
\end{align}
is the loss function for the \textbf{Forward problems}.

Unconstrained minimization methods are usually terminated when the gradient of the objective function is  small but not necessarily zero. Set the termination criterion for iterations as
\begin{equation}
    |\nabla_{\bm{\theta}} \mathcal{L}(\bm{\theta},\bm{\lambda}; \mu)| \leq \delta,
\end{equation}
where $\delta$ is a given scalar representing a small number.

\begin{figure}[htbp]
    \centering
    \includegraphics[width=0.7\textwidth]{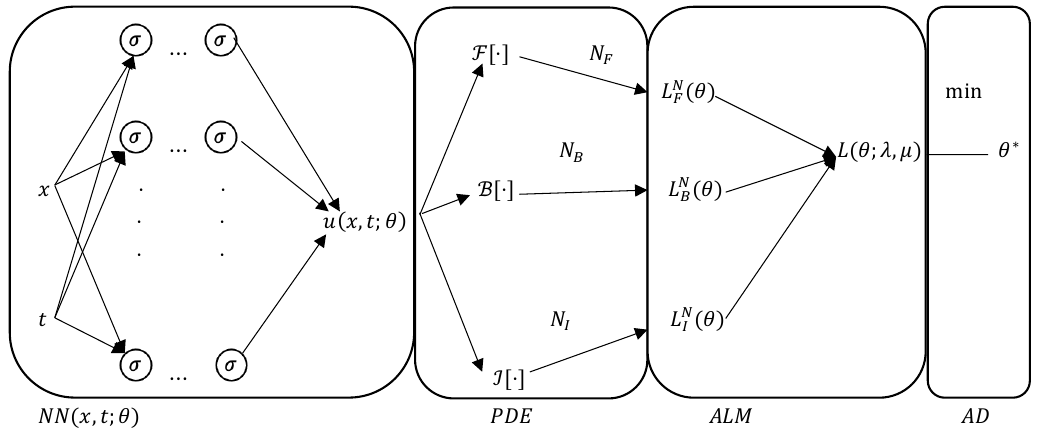}
    \caption{The schematic diagram of the ALM-PINNs program for solving forward problems is presented. NN denotes the neural network module, PDE represents the module for incorporating the physical model, ALM indicates the module for constructing the loss function, and AD refers to the automatic differentiation module.}
    \label{fig:11}
\end{figure}

We summarize the above contents in Algorithm \ref{alg1:Framework}, along with the schematic diagram \ref{fig:11} of the ALM-PINNs program for solving forward problems.

\begin{algorithm}[htbp]  
	\caption{ALM-PINNs Iterative Procedure for Solving PDEs}  
	\label{alg1:Framework}
	\SetAlgoLined
	\KwIn { The set of governing equation points \(X_r\), boundary and initial data points \(X_{\text{data}}\), observed data \(u_{\text{data}}\), boundary condition weight $\lambda_{\mathcal{B}}$, initial condition weight $\lambda_{\mathcal{I}}$, penalty coefficient \(\mu\), Initial neural network parameters: $\bm{\theta}_0$;}
	\KwOut {Updated model parameters \( \bm{\theta^* }\), loss history}
 
1. Initialize model \( \bm{\theta} \) \;
2. Set learning rate schedule \( lr \) as PiecewiseConstantDecay \;
3. Define optimizer \( optim \) as Adam with learning rate \( lr \) \;
4. Function ComputeLoss(\( \bm{\theta}, X_r, X_{\text{data}}, u_{\text{data}}, \lambda_{\mathcal{B}}, \lambda_{\mathcal{I}}, \mu \)):
The total loss function is calculated using Eq. (\ref{eq:ALM_Forward}). \;
5. Function ComputeGradient(\( \bm{\theta}, X_r, X_{\text{data}}, u_{\text{data}}, \lambda_{\mathcal{B}}, \lambda_{\mathcal{I}}, \mu \)):
    The loss and related terms are computed using the ComputeLoss function, while the gradients of the loss with respect to the parameters $\bm{\theta}$ are calculated via automatic differentiation in TensorFlow. \;
6. Function TrainStep():
    Compute loss and gradients using ComputeGradient. Update model parameters \( \theta \) using optimizer \( Adam \) \;
    \For{$i \leftarrow 1$ \KwTo Epoch}{
        \For{$j \leftarrow 1$ \KwTo Batch}{
		 $\text{loss}, \text{penalty}, \text{gover\_loss}, \text{bc\_ls}, \text{ic\_ls} \leftarrow \text{TrainStep()}$\;
        \If {$\sqrt{\text{penalty}} > \epsilon$}{
            Update the penalty parameter: $\mu \leftarrow \min(\text{penalty} \cdot \mu, \mu_{\text{max}})$\;
            Update boundary Lagrange multiplier: $\lambda_{\mathcal{B}} \leftarrow \lambda_{\mathcal{B}} + \mu \cdot \text{bc\_ls}$ \;
            Update initial condition Lagrange multiplier: $\lambda_{\mathcal{I}} \leftarrow \lambda_{\mathcal{I}} + \mu \cdot \text{ic\_ls}$ \;
            Append the current loss to history: $\text{hist.append}(\text{loss})$ \;
        }
        
        \If{$\text{gover\_loss} < \text{best\_gover\_loss}$}{
            Update best loss: $\text{best\_loss} \leftarrow \text{loss}$ \;
            Update best governing equation loss: $\text{best\_gover\_loss} \leftarrow \text{gover\_loss}$ \;
            Store the current model and loss history: $\text{best\_model} \leftarrow \text{current model}$
            $\text{best\_hist} \leftarrow \text{hist.copy()}$ \;
            Record the best computation time: $\text{best\_time} \leftarrow \text{time()} - t_0$
        }
		}
    }
\end{algorithm}

\section{The construction principle of the ALM-PINNs algorithm for solving inverse problems}
\label{sec:ALM-PINNs-inverse}

\subsection{Problem setup}

\textbf{The parameter inversion problem discussed in this section is as follows.} For the differential equation defined in (\ref{gover_Eq}), assuming that the parameter $\bm{v}$ is unknown, we aim to infer the unknown parameters by constructing an additional dataset $T_j\in [0,t_T], j=1,\dots,n, \{(\bm{x}^{(i)}_{inner}, T_j), u^{(i)}\}_{i=1}^{i=N_{\mathcal{E}}}$. The smaller the value of $n$, the better. This dataset consists of the exact values of $u$ at arbitrary measurement points $(\bm{x},t)$ within the entire domain or values perturbed by noise following different distributions (see Equation (\ref{perturbed_u})).
\begin{equation}
    \label{perturbed_u}
    \tilde{u}(\bm{x}^{(i)},t^{(i)}) = u(\bm{x}^{(i)},t^{(i)};\bm{\theta},\bm{v}) + \epsilon^{(i)} ,\quad i = 1,\cdots, N_{\mathcal{E}},
\end{equation}
where $\bm{\epsilon}$ is a random variable that follows a normal distribution, Laplace distribution, and log-normal distribution respectively. The theorem in the next subsection will illustrate the construction rule for the objective functional on the additional dataset when the data follows different prior distributions.

\begin{figure}[htbp]
    \centering
    \includegraphics[width=0.5\linewidth]{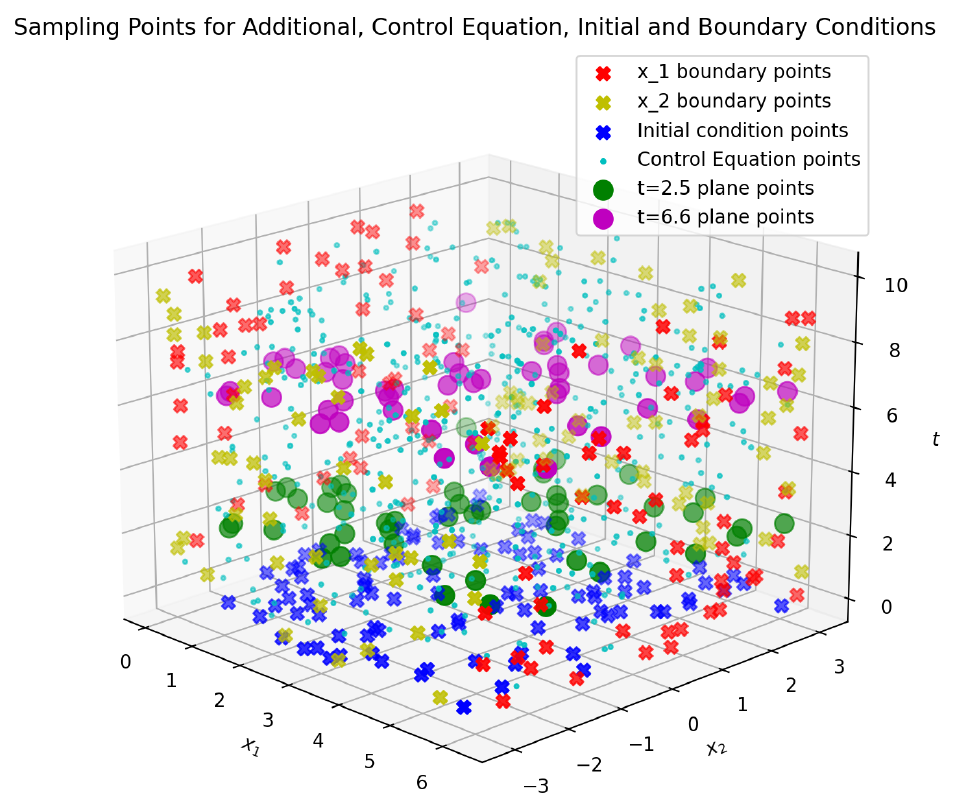}
    \caption{Sampling Points for Additional, Control Equation, Initial and Boundary Conditions}
    \label{fig:additional}
\end{figure}

Figure \ref{fig:additional} shows the sampling diagram of additional data points at times 
$t=2.5$ and $t=6.6$, along with the selected points within the control equation and the initial/boundary conditions. These data will be used in the inverse problem to identify the parameter $\bm{v}$ and the neural network parameters $\bm{\theta}$.

\subsection{The construction rule for the objective functional on the additional dataset}

\begin{theorem}
\label{thm}
    When solving inverse problems using ALM-PINNs algorithm, assuming additional noisy data set $\bm{\varepsilon}$, $\tilde{u}(\bm{x}^{(i)}, t^{(i)}) = u(\bm{x}^{(i)}, t^{(i)};\bm{\theta}, \bm{v}) + \varepsilon^{(i)} \quad  i = 1,\cdots, N_{\mathcal{M}}$, $\tilde{u}(\bm{x}^{(i)}, t^{(i)})$ is the observed data at $(\bm{x}^{(i)}, t^{(i)})$, $u(\bm{x}^{(i)}, t^{(i)};\bm{\theta}, \bm{v})$ is the predicted data of the neural network function at $(\bm{x}^{(i)}, t^{(i)})$, $\varepsilon^{(i)}$ is the error data associated with $(\bm{x}^{(i)}, t^{(i)})$, and $N_{\mathcal{M}}$ is the total number of observations. Then, the construction rule of the additional data term $\mathcal{L}_{\mathcal{E}}^{\mathcal{N}}(\bm{\theta}, \bm{v})$ in the loss function can be described as follows

    \rm (1) If $\bm{\varepsilon} \sim N(\bm{0}, \sigma^2\bm{I})$, then we choose 
    \begin{equation*}
        \mathcal{L}_{\mathcal{E}}^{\mathcal{N}}({\bm{\theta}, \bm{v}}) =  \frac{1}{N_{\mathcal{M}}}\sum_{i=1}^{N_{\mathcal{M}}}\frac{(\tilde{u}(\bm{x}^{(i)}, t^{(i)}) - u(\bm{x}^{(i)}, t^{(i)};\bm{\theta}, \bm{v}))^2}{2\sigma^2}.
    \end{equation*}
    Specifically, if $\bm{\varepsilon} \sim N(\bm{0}, \frac{1}{2}\bm{I})$, $\mathcal{L}_{\mathcal{E}}^{\mathcal{N}}({\bm{\theta}, \bm{v}})$ takes the $2$-norm form.

    \rm (2) If $\bm{\varepsilon} \sim Laplace(\bm{0}, 2\gamma^2\bm{I})$, then we choose 
    \begin{equation*}
        \mathcal{L}_{\mathcal{E}}^{\mathcal{N}}({\bm{\theta}, \bm{v}}) =  \frac{1}{N_{\mathcal{M}}}\sum_{i=1}^{N_{\mathcal{M}}}\frac{|\tilde{u}(\bm{x}^{(i)}, t^{(i)}) - u(\bm{x}^{(i)}, t^{(i)};\bm{\theta}, \bm{v})|}{\gamma}.
    \end{equation*}
    Specifically, if $\bm{\varepsilon} \sim Laplace(\bm{0}, 2\bm{I})$, $\mathcal{L}_{\mathcal{E}}^{\mathcal{N}}({\bm{\theta}, \bm{v}})$ takes the $1$-norm form.

    \rm (3) If $\bm{\varepsilon} \sim LogN(\bm{0}, \sigma^2\bm{I})$, then we choose 
   \begin{equation*}
    \begin{aligned}
        \mathcal{L}_{\mathcal{E}}^{\mathcal{N}}(\bm{\theta}, \bm{v}) = & \frac{1}{N_{\mathcal{M}}} \sum_{i=1}^{N_{\mathcal{M}}} \left[ \ln\left(\sqrt{2\pi}\sigma \left( \tilde{u}(\bm{x}^{(i)}, t^{(i)}) - u(\bm{x}^{(i)}, t^{(i)};\bm{\theta}, \bm{v})\right)\right) \right. \\
        & \left. + \frac{\left( \ln\left(\tilde{u}(\bm{x}^{(i)}, t^{(i)}) - u(\bm{x}^{(i)}, t^{(i)};\bm{\theta}, \bm{v})\right) \right)^2}{2\sigma^2} \right].
    \end{aligned}
\end{equation*}
    \label{Thm:3.1}
\end{theorem}

\begin{remark}
    According to the derivation in this paper for three scenarios, regardless of the distribution of prior data errors, the construction form of the objective functional on the dataset can be derived once the mean and variance are known.
\end{remark}

\begin{proof}

\begin{case} 
Assume the error $\bm{\varepsilon}$ is normally distributed with mean 0 and standard deviation $\sigma$.
\end{case}
The probability density for $d_i$ takes the form of 
\begin{equation}
    f_i(d_i|\bm{v})=f(\varepsilon^{(i)}) = \frac{1}{\sigma_i \sqrt{2\pi}}\exp(-\frac{\varepsilon^{(i)}}{2\sigma_i^2}). 
\end{equation}
The likelihood function for the complete data set is the product of the individual likelihoods
\begin{equation}
L(\bm{\theta},\bm{v}|\bm{d})=\prod_{i=1}^{N_{\mathcal{M}}}\frac{1}{\sigma_i \sqrt{2\pi}}\exp(-\frac{(\tilde{u}(\bm{x}^{(i)}, t^{(i)}) - u(\bm{x}^{(i)}, t^{(i)});\bm{\theta}, \bm{v})^2}{2\sigma_i^2}).
\end{equation}
The natural logarithm likelihood function is 
\begin{equation}
    \begin{aligned}
        l(\bm{\theta}, \bm{v}|\bm{d}) &= \ln(\prod_{i=1}^{N_{\mathcal{M}}}\frac{1}{\sigma_i \sqrt{2\pi}}\exp(-\frac{(\tilde{u}(\bm{x}^{(i)}, t^{(i)}) - u(\bm{x}^{(i)}, t^{(i)});\bm{\theta}, \bm{v})^2}{2\sigma_i^2})) \\
            &= \sum_{i=1}^{N_{\mathcal{M}}}\ln(\frac{1}{\sigma_i \sqrt{2\pi}}\exp(-\frac{(\tilde{u}(\bm{x}^{(i)}, t^{(i)}) - u(\bm{x}^{(i)}, t^{(i)};\bm{\theta}, \bm{v}))^2}{2\sigma_i^2})) \\
            &=  \sum_{i=1}^{N_{\mathcal{M}}} \ln(\frac{1}{\sigma_i \sqrt{2\pi}})-\sum_{i=1}^{N_{\mathcal{M}}}\frac{(\tilde{u}(\bm{x}^{(i)}, t^{(i)}) - u(\bm{x}^{(i)}, t^{(i)};\bm{\theta}, \bm{v}))^2}{2\sigma_i^2}. \\
\end{aligned}
\end{equation}
By ignore the constant factor of $\frac{1}{2}$, the problem becomes
\begin{equation}
\begin{aligned}
\label{3}
    \max_{\bm{\theta}, \bm{v}|\bm{d}}  l(\bm{\theta}, \bm{v})&\Rightarrow \max_{\bm{\theta}, \bm{v}|\bm{d}}  -\sum_{i=1}^{N_{\mathcal{M}}}\frac{(\tilde{u}(\bm{x}^{(i)}, t^{(i)}) - u(\bm{x}^{(i)}, t^{(i)};\bm{\theta}, \bm{v}))^2}{\sigma_i^2} \\
    &\Rightarrow \min_{\bm{\theta}, \bm{v}|\bm{d}}  \sum_{i=1}^{N_{\mathcal{M}}}\frac{(\tilde{u}(\bm{x}^{(i)}, t^{(i)}) - u(\bm{x}^{(i)}, t^{(i)};\bm{\theta}, \bm{v}))^2}{\sigma_i^2},
\end{aligned}
\end{equation} 
which is commonly referred to as the weighted least squares problem.

\begin{case}
    Assuming the error $\bm{\varepsilon}$ follows a Laplace distribution with mean $0$ and variance $2\gamma^2$.
\end{case}
Expressing $\bm{\varepsilon}$ as
\begin{equation}
    \bm{\varepsilon} \sim Laplace(\bm{0}, 2\gamma^2\bm{I}),
\end{equation}
and its probability density function is
\begin{equation}
    f(\varepsilon^{(i)}) = \frac{1}{2\gamma}\exp\left(-\frac{|\varepsilon^{(i)}|}{\gamma}\right).
\end{equation}
Thus the conditional probability density function of $\tilde{u}$ is
\begin{equation}
    f(\tilde{u}|x^{(i)}, t^{(i)};\bm{\theta}, \bm{v})=\frac{1}{2\gamma}\exp\left(-\frac{|\tilde{u}(\bm{x}^{(i)}, t^{(i)}) - u(\bm{x}^{(i)}, t^{(i)};\bm{\theta}, \bm{v})|}{\gamma}\right).
\end{equation}
According to the maximum likelihood estimation principle, its likelihood function is
\begin{equation}
    L(\bm{\theta}, \bm{v})=\prod_{i=1}^{N_{\mathcal{M}}}\frac{1}{2\gamma}\exp\left(-\frac{|\tilde{u}(\bm{x}^{(i)}, t^{(i)}) - u(\bm{x}^{(i)}, t^{(i)};\bm{\theta}, \bm{v})|}{\gamma}\right),
\end{equation}
and its log-likelihood function is
\begin{equation}
    \begin{aligned}
        l(\bm{\theta}, \bm{v}) &= \ln\left(\prod_{i=1}^{N_{\mathcal{M}}}\frac{1}{2\gamma}\exp\left(-\frac{|\tilde{u}(\bm{x}^{(i)}, t^{(i)}) - u(\bm{x}^{(i)}, t^{(i)};\bm{\theta}, \bm{v})|}{\gamma}\right)\right) \\
            &= \sum_{i=1}^{N_{\mathcal{M}}}\ln\left(\frac{1}{2\gamma}\exp\left(-\frac{|\tilde{u}(\bm{x}^{(i)}, t^{(i)}) - u(\bm{x}^{(i)}, t^{(i)};\bm{\theta}, \bm{v})|}{\gamma}\right)\right) \\
            &=  N_{\mathcal{M}} \ln\left(\frac{1}{2\gamma}\right)-\sum_{i=1}^{N_{\mathcal{M}}}\frac{|\tilde{u}(\bm{x}^{(i)}, t^{(i)}) - u(\bm{x}^{(i)}, t^{(i)};\bm{\theta}, \bm{v})|}{\gamma}.
    \end{aligned}
\end{equation}
Therefore,
\begin{equation}
    \begin{aligned}
        \max_{\bm{\theta}, \bm{v}}  l(\bm{\theta}, \bm{v})&\Rightarrow \max_{\bm{\theta}, \bm{v}}  -\sum_{i=1}^{N_{\mathcal{M}}}\frac{|\tilde{u}(\bm{x}^{(i)}, t^{(i)}) - u(\bm{x}^{(i)}, t^{(i)};\bm{\theta}, \bm{v})|}{\gamma} \\
        &\Rightarrow \min_{\bm{\theta}, \bm{v}}  \sum_{i=1}^{N_{\mathcal{M}}}\frac{|\tilde{u}(\bm{x}^{(i)}, t^{(i)}) - u(\bm{x}^{(i)}, t^{(i)};\bm{\theta}, \bm{v})|}{\gamma}.
    \end{aligned}
\end{equation}
At this point, let
\begin{equation}
    \mathcal{L}_{\mathcal{E}}^{\mathcal{N}}({\bm{\theta}, \bm{v}}) =  \frac{1}{N_{\mathcal{M}}}\sum_{i=1}^{N_{\mathcal{M}}}\frac{|\tilde{u}(\bm{x}^{(i)}, t^{(i)}) - u(\bm{x}^{(i)}, t^{(i)};\bm{\theta}, \bm{v})|}{\gamma}.
    \label{3. 47}
\end{equation}

Noting that if $\gamma = 1$, then
\begin{equation}
    \bm{\varepsilon} \sim Laplace(\bm{0}, 2\bm{I}).
\end{equation}
Equation (\ref{3. 47}) takes the $1$-norm form, which indirectly confirms that when designing the objective functional, the form of the data terms should be chosen based on prior error estimates for different errors.

\begin{case}
    Assuming the error $\bm{\varepsilon}$ follows a log-normal distribution with mean $0$ and variance $\sigma^2$.
\end{case}
Expressing $\bm{\varepsilon}$ as
\begin{equation}
    \bm{\varepsilon} \sim LogN(\bm{0}, \sigma^2\bm{I}),
\end{equation}
we know the skewed distribution probability density function is
\begin{equation}
    f(\varepsilon^{(i)}) = \left\{
    \begin{array}{lll}
        \frac{1}{\sqrt{2\pi}\sigma\varepsilon^{(i)}}\exp\left(-\frac{(\ln(\varepsilon^{(i)}))^2}{2\sigma^2}\right)   &      & \varepsilon^{(i)} >0\\
        0&      & \varepsilon^{(i)}  \leq 0
    \end{array} \right.  ,
\end{equation}
thus the conditional probability density function of $\tilde{u}$ is
\begin{equation}
    \begin{aligned}
        &f(\tilde{u}|x^{(i)}, t^{(i)};\bm{\theta}, \bm{v}) \\
        &=\frac{1}{\sqrt{2\pi}\sigma(\tilde{u}(\bm{x}^{(i)}, t^{(i)}) - u(\bm{x}^{(i)}, t^{(i)};\bm{\theta}, \bm{v}))}\exp\left(-\frac{(\ln(\tilde{u}(\bm{x}^{(i)}, t^{(i)}) - u(\bm{x}^{(i)}, t^{(i)};\bm{\theta}, \bm{v})))^2}{2\sigma^2}\right).
    \end{aligned}
\end{equation}
According to the maximum likelihood estimation principle, its likelihood function is
\begin{equation}
    \begin{aligned}
        &L(\bm{\theta}, \bm{v})\\
        &=\prod_{i=1}^{N_{\mathcal{M}}}\frac{1}{\sqrt{2\pi}\sigma(\tilde{u}(\bm{x}^{(i)}, t^{(i)}) - u(\bm{x}^{(i)}, t^{(i)};\bm{\theta}, \bm{v}))}\exp\left(-\frac{(\ln(\tilde{u}(\bm{x}^{(i)}, t^{(i)}) - u(\bm{x}^{(i)}, t^{(i)};\bm{\theta}, \bm{v})))^2}{2\sigma^2}\right),
    \end{aligned}
\end{equation}
and its log-likelihood function is
\begin{equation}
    \begin{aligned}
        &l(\bm{\theta}, \bm{v}) \\
            &= \ln\left(\prod_{i=1}^{N_{\mathcal{M}}}\frac{1}{\sqrt{2\pi}\sigma(\tilde{u}(\bm{x}^{(i)}, t^{(i)}) - u(\bm{x}^{(i)}, t^{(i)};\bm{\theta}, \bm{v}))}\exp\left(-\frac{(\ln(\tilde{u}(\bm{x}^{(i)}, t^{(i)}) - u(\bm{x}^{(i)}, t^{(i)};\bm{\theta}, \bm{v})))^2}{2\sigma^2}\right)\right)\\
            &= \sum_{i=1}^{N_{\mathcal{M}}}\ln\left(\frac{1}{\sqrt{2\pi}\sigma(\tilde{u}(\bm{x}^{(i)}, t^{(i)}) - u(\bm{x}^{(i)}, t^{(i)};\bm{\theta}, \bm{v}))}\exp\left(-\frac{(\ln(\tilde{u}(\bm{x}^{(i)}, t^{(i)}) - u(\bm{x}^{(i)}, t^{(i)};\bm{\theta}, \bm{v})))^2}{2\sigma^2}\right)\right)\\
            &= \sum_{i=1}^{N_{\mathcal{M}}}-\ln\left(\sqrt{2\pi}\sigma(\tilde{u}(\bm{x}^{(i)}, t^{(i)}) - u(\bm{x}^{(i)}, t^{(i)};\bm{\theta}, \bm{v}))\right)-\frac{(\ln(\tilde{u}(\bm{x}^{(i)}, t^{(i)})-u(\bm{x}^{(i)}, t^{(i)};\bm{\theta}, \bm{v})))^2}{2\sigma^2}. \\
    \end{aligned}
\end{equation}
Therefore,
\begin{equation}
    \begin{aligned}
        &\max_{\bm{\theta}, \bm{v}}  l(\bm{\theta}, \bm{v}) \Rightarrow\\
        &\min_{\bm{\theta}, \bm{v}}  \sum_{i=1}^{N_{\mathcal{M}}}\ln\left(\sqrt{2\pi}\sigma(\tilde{u}(\bm{x}^{(i)}, t^{(i)}) - u(\bm{x}^{(i)}, t^{(i)};\bm{\theta}, \bm{v}))\right)+\frac{(\ln(\tilde{u}(\bm{x}^{(i)}, t^{(i)})-u(\bm{x}^{(i)}, t^{(i)};\bm{\theta}, \bm{v})))^2}{2\sigma^2}. 
    \end{aligned}
\end{equation}
At this point, let
\begin{equation}
    \begin{aligned}
        &\mathcal{L}_{\mathcal{E}}^{\mathcal{N}}(\bm{\theta}, \bm{v}) =\\
        &\frac{1}{N_{\mathcal{M}}}\sum_{i=1}^{N_{\mathcal{M}}}\ln\left(\sqrt{2\pi}\sigma(\tilde{u}(\bm{x}^{(i)}, t^{(i)}) - u(\bm{x}^{(i)}, t^{(i)};\bm{\theta}, \bm{v}))\right)+\frac{(\ln(\tilde{u}(\bm{x}^{(i)}, t^{(i)})-u(\bm{x}^{(i)}, t^{(i)};\bm{\theta}, \bm{v})))^2}{2\sigma^2}. 
    \end{aligned}
\end{equation}
\end{proof}


Unlike the optimization problem established in Problem \ref{eq:primal}, this paper formulates the boundary term and initial condition term of the inverse problem as the minimization objective, while treating the control equation term and data term as constraint terms. The rationale behind this is that, during the neural network iterations, we are more concerned with the equation parameters \(v_1\) and \(v_2\) rather than the neural network parameters themselves. In accordance with the construction rules of the augmented Lagrangian method (ALM), weights are applied to \(v_1\) and \(v_2\) to accelerate their updates during the iteration process.

\begin{equation}
\label{eq:inverse}
    \begin{aligned} 
    &\min_{\bm{\theta}\in \mathbb{R}^n} \quad \mathcal{L}_{\mathcal{B}}^{\mathcal{N}}(\bm{\theta})+\mathcal{L}_{\mathcal{I}}^{\mathcal{N}}(\bm{\theta}),   \\
    &\begin{array}{r@{\quad}r@{}l@{\quad}l}
    s. t.   &    \mathcal{L}_{\mathcal{F}}^{\mathcal{N}}(\bm{\theta},\bm{v})=0, \\
    &   \mathcal{L}_{\mathcal{E}}^{\mathcal{N}}(\bm{\theta},\bm{v})=0, \\
    \end{array}
    \end{aligned}
\end{equation}

Based on Section \ref{sec:ALM-PINNs-forward} and the analysis in this section, the objective function for the \textbf{inverse problems} is defined as:
\begin{align}
    \label{Eq:ALM_inverse}
    \mathcal{L}_{Inverse}(\bm{\theta},\bm{v};\bm{\lambda}, \bm{\mu})&=\mathcal{L}_{\mathcal{B}}^{\mathcal{N}}(\bm{\theta})+\mathcal{L}_{\mathcal{I}}^{\mathcal{N}}(\bm{\theta})+\bm{\lambda}_{\mathcal{F}}\mathcal{L}_{\mathcal{F}}^{\mathcal{N}}(\bm{\theta},\bm{v})+\bm{\lambda}_{\mathcal{E}}\mathcal{L}_{\mathcal{E}}^{\mathcal{N}}(\bm{\theta},\bm{v})\\
    &+\frac{\mu_{\mathcal{F}}}{2}|\mathcal{L}_{\mathcal{F}}^{\mathcal{N}}(\bm{\theta},\bm{v})|^2+\frac{\mu_{\mathcal{E}}}{2}|\mathcal{L}_{\mathcal{E}}^{\mathcal{N}}(\bm{\theta},\bm{v})|^2, 
\end{align}

\subsection{Strategies for overcoming the ill-posedness in Neural Networks Training for Inverse Problems}
\label{ill-posed}

Inverse problems are typically ill-posed, and current research focuses on both the numerical reproduction of solutions and the theoretical proof of model uniqueness and stability. One of the key reasons for the popularity of Physics-Informed Neural Networks (PINNs) is that they handle forward and inverse problems with comparable difficulty. Specifically, using the loss function model established for the forward problem, solving the inverse problem only requires introducing the parameters to be estimated and iteratively training them along with the neural network’s weights and biases. However, this approach presents some inherent challenges: 

First, neural networks (NNs) have a large number of parameters, and it is difficult to prioritize the training of the target parameters in the iterative process. Second, the training process of NNs cannot alter the inherent ill-posedness of the model, which often leads to lower accuracy in solving inverse problems compared to forward problems.

To address these issues, this paper introduces prior information on the solution from four perspectives to improve the accuracy of inverse problem solutions:

\textbf{(1). Introducing a prior distribution for the error}: By constructing a loss function for the inverse problem using the maximum likelihood estimation principle combined with the ALM-PINNs algorithm, different prior information about the error is incorporated into different loss functions to improve the accuracy of the inferred parameters.

\textbf{(2). Introducing constraints on the solution range}: By imposing certain prior restrictions on the solution variables, the parameters are trained within a specific range. For example, human height must be positive, and the diameter of a mature sycamore tree falls within a certain range. This approach can accelerate model convergence and improve solution accuracy.

\textbf{(3). Incorporating a pre-trained model from the forward problem}: Neural networks can be used to solve differential equations, and forward and inverse problems are often two sides of the same coin. The pre-trained model (i.e., the weights and biases under a specific neural network structure) from the forward problem, saved as an h5 file, can be introduced into the inverse problem model. Furthermore, the accuracy of the forward problem solution can influence the quality of the inferred parameters in the inverse problem.

\textbf{(4). Weighting the parameters of interest}: In this approach, the neural network’s weight parameters are only fine-tuned, and the focus of the training is shifted toward the parameters of interest. The original PINNs method treats the target parameters as exogenous variables in the training process, training them along with the neural network. For example, if the model itself requires training 12,921 parameters and there are 2 additional exogenous variables, the total number of parameters for the inverse problem would be 12,923, although we are only interested in 2 of them. Therefore, during the inverse problem-solving process, we apply weights to the parameters of interest and aim for an optimal solution within their allowable range.

The numerical experiments for the inverse problem in this paper are designed based on the above four aspects. Specifically, the pre-trained model (h5 file) used in the ALM-PINNs algorithm for the inverse problem is obtained from solving the forward problem using ALM-PINNs. Similarly, the pre-trained model (h5 file) used in the PINNs algorithm for comparison is obtained from solving the forward problem using PINNs.

The schematic diagrams of the ALM-PINNs algorithm for solving inverse problems and ALM-PINNs Iterative Procedure for IPs are shown in Figure \ref{fig:11} and Algorithm \ref{alg2:Framework}, respectively. 

\begin{algorithm}[htbp]  
	\caption{ALM-PINNs Iterative Procedure for IPs}  
	\label{alg2:Framework}
	\SetAlgoLined
	\KwIn { The set of governing equation points \(X_r\), boundary and initial data points \(X_{\text{data}}\), observed data \(u_{\text{data}}\), boundary condition weight $\lambda_{\mathcal{B}}$, initial condition weight $\lambda_{\mathcal{I}}$, penalty coefficient \(\mu\), Pre-trained model $\text{pretrained\_model}$ $\bm{\theta}^*$;}
	\KwOut {Updated model parameters \( \hat{\bm{\theta }}, \bm{v}\)}
 
1. Load the pre-trained model $\bm{\theta^*}$ and add new hidden layers as appropriate\;
2. Set learning rate schedule \( lr \) as PiecewiseConstantDecay \;
3. Define optimizer \( optim \) as Adam with learning rate \( lr \) \;
4. Initialize parameters $\bm{v}$ with a starting value and use tf.clip\_by\_value to ensure it stays within predefined bounds during optimization. \;
5. Function ComputeLoss(\( \bm{\theta}^*, X_r, X_{\text{data}}, u_{\text{data}, \bm{v}}, \lambda_{\mathcal{B}}, \lambda_{\mathcal{I}}, \mu \)): 
The total loss function is calculated using Eq. (\ref{Eq:ALM_inverse}). \;
6. Function ComputeGradient(\( \bm{\theta}^*, X_r, X_{\text{data}}, u_{\text{data}}, \bm{v}, \lambda_{\mathcal{B}}, \lambda_{\mathcal{I}}, \mu \)):
    The loss and related terms are computed using the ComputeLoss function, while the gradients of the loss with respect to the parameters $\bm{\theta}$ and $\bm{v}$ are calculated via automatic differentiation in TensorFlow. \;
7. Function TrainStep():
    Compute loss and gradients using ComputeGradient. Update model parameters \( \bm{\theta},\bm{v} \) using optimizer \(Adam \) \;
    \For{$i \leftarrow 1$ \KwTo Epoch}{
        \For{$j \leftarrow 1$ \KwTo Batch}{
		$\text{loss}, \text{penalty}, \text{gover\_loss}, \text{bc\_ls}, \text{ic\_ls} \leftarrow \text{TrainStep()}$\;
        \If {$\sqrt{\text{penalty}} > \epsilon$}{
            Update the penalty parameter: $\mu \leftarrow \min(\text{penalty} \cdot \mu, \mu_{\text{max}})$\;
            Update boundary Lagrange multiplier: $\lambda_{\mathcal{B}} \leftarrow \lambda_{\mathcal{B}} + \mu \cdot \text{bc\_ls}$ \;
            Update initial condition Lagrange multiplier: $\lambda_{\mathcal{I}} \leftarrow \lambda_{\mathcal{I}} + \mu \cdot \text{ic\_ls}$ \;
            Append the current loss to history: $\text{hist.append}(\text{loss})$ \;
        }
        
        \If{$\text{gover\_loss} < \text{best\_gover\_loss}$}{
            Update best loss: $\text{best\_loss} \leftarrow \text{loss}$ \;
            Update best governing equation loss: $\text{best\_gover\_loss} \leftarrow \text{gover\_loss}$ \;
            Store the current model and loss history: $\text{best\_model} \leftarrow \text{current model}$
            $\text{best\_hist} \leftarrow \text{hist.copy()}$ \;
            Record the best computation time: $\text{best\_time} \leftarrow \text{time()} - t_0$
        }
    }
}
\end{algorithm}

\begin{figure}[htbp]
    \centering
    \includegraphics[width=0.7\textwidth]{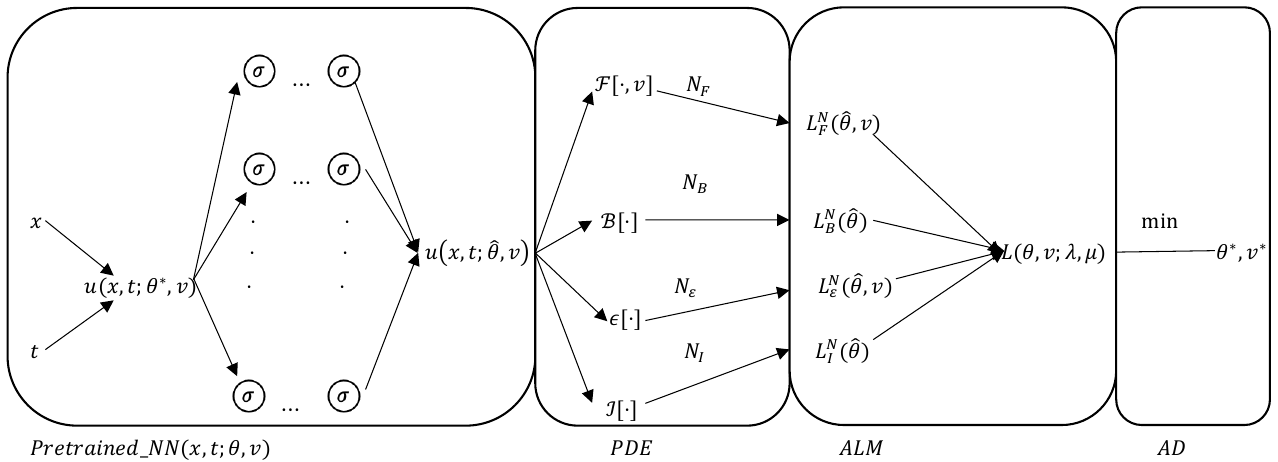}
    \caption{The schematic diagram of the ALM-PINNs program for solving inverse problems is presented. $Pretrained_{NN}$ denotes the neural network module incorporating the pretrained model, PDE represents the module for introducing the physical model, ALM indicates the module for constructing the loss function, and AD refers to the automatic differentiation module.}
    \label{fig:22}
\end{figure}



\section{Numerical examples of solving nonlinear differential equations}
\label{sec:forward}
In this section, we will apply both the standard PINNs and the ALM-PINNs methods to solve specific partial differential equations with known exact solutions, aiming to demonstrate that our approach achieves more accurate numerical results. This study evaluates the accuracy of the models by calculating the global relative error $\epsilon_r(\hat{u},u)$, maximum absolute error $\epsilon_\infty(\hat{u},u)$, and mean absolute error $\epsilon_\infty(\hat{u},u)$ of the numerical solutions obtained from PINNs and ALM-PINNs respectively, 
\begin{equation*}
\epsilon_r(\hat{u},u) = \frac{\|\hat{u} - u\|_2}{\|u\|_2},
\quad
\epsilon_\infty(\hat{u},u)= \max |\hat{u} - u|,
\quad
\epsilon_a(\hat{u},u) = \frac{\sum |\hat{u} - u|}{n},
\end{equation*}
where $n$ represents the total number of data points.

 \subsection{The 1+1 dimensional nonlinear equation}
\label{eg1}
The first type of nonlinear equations has important applications in various fields, including ecology, chemical reaction dynamics, and fluid mechanics. For instance, it can be used to model the diffusion and competition of biological populations or to describe the transfer of heat in different media.
\begin{equation}
    \left\{
    \begin{array}{lll}
    u_t=2(u_x)^2+2uu_{xx}+u-u^2,      &      & t\in [0, 4], x \in [0, 1],\\
    u(x, 0)=\left\{\frac{\tanh(-\frac{x}{4})}{2}+\frac{1}{2}\right\}^{-1},&      & x \in [0, 1],\\
    u(0, t)=\left\{\frac{\tanh(\frac{t}{4})}{2}+\frac{1}{2}\right\}^{-1},   &      & t\in [0, 4],\\
    u(1, t)=\left\{\frac{\tanh(\frac{t}{4}-\frac{1}{4})}{2}+\frac{1}{2}\right\}^{-1},   &      & t\in [0, 4].
    \end{array} \right.  
\end{equation}
And the equation has an exact solution
\begin{equation}
    u(x, t)= \left\{\frac{1}{2}+\frac{\tanh(\frac{t}{4}-\frac{x}{4})}{2}\right\}^{-1}.
\end{equation}
The network architecture selected for this example is detailed in Table \ref{table:compare_1}. It features an input layer with 2 nodes and an output layer with 1 node. Additionally, there is 1 scaling layer and 7 hidden layers, each containing 40 nodes, with the activation function configured to the hyperbolic tangent (tanh) function. The learning rate is defined in a segmented manner as follows: for iterations in the range \([0, 100]\), the learning rate is set to \(1 \times 10^{-2}\); for iterations in the range \([100, 1000]\), it is \(1 \times 10^{-3}\); for iterations in the range \([1000, 2500]\), it is \(5 \times 10^{-4}\); and for iterations beyond \(2500\), the learning rate is adjusted to \(1 \times 10^{-4}\).

\begin{figure}[htbp]
    \centering
    \begin{subfigure}{0.3\textwidth}
        \centering
        \includegraphics[width=\textwidth]{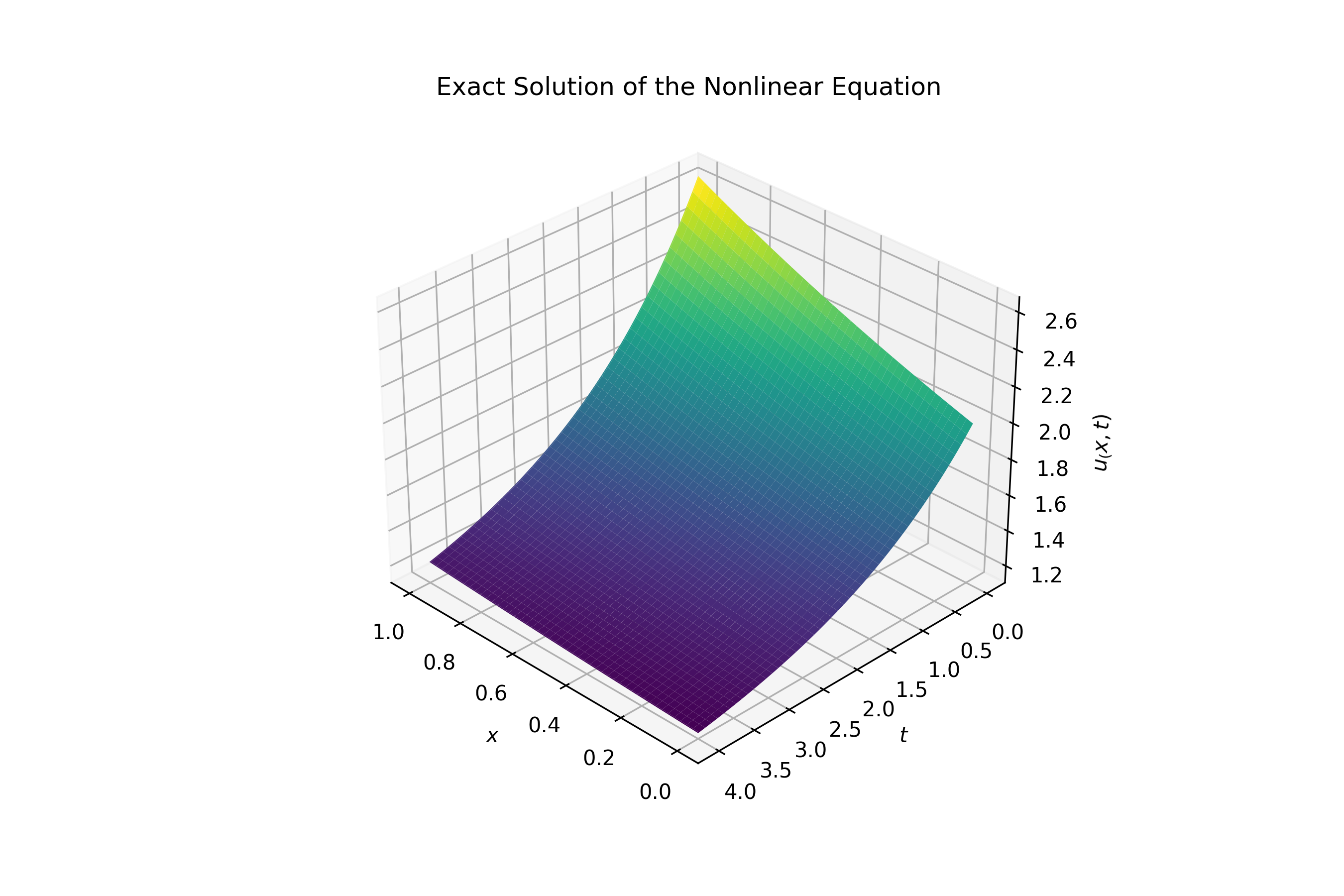}
        \caption{Exact solution for the nonlinear equation}
        \label{fig:sub2}
    \end{subfigure}
    \begin{subfigure}{0.3\textwidth}
        \centering
        \includegraphics[width=\textwidth]{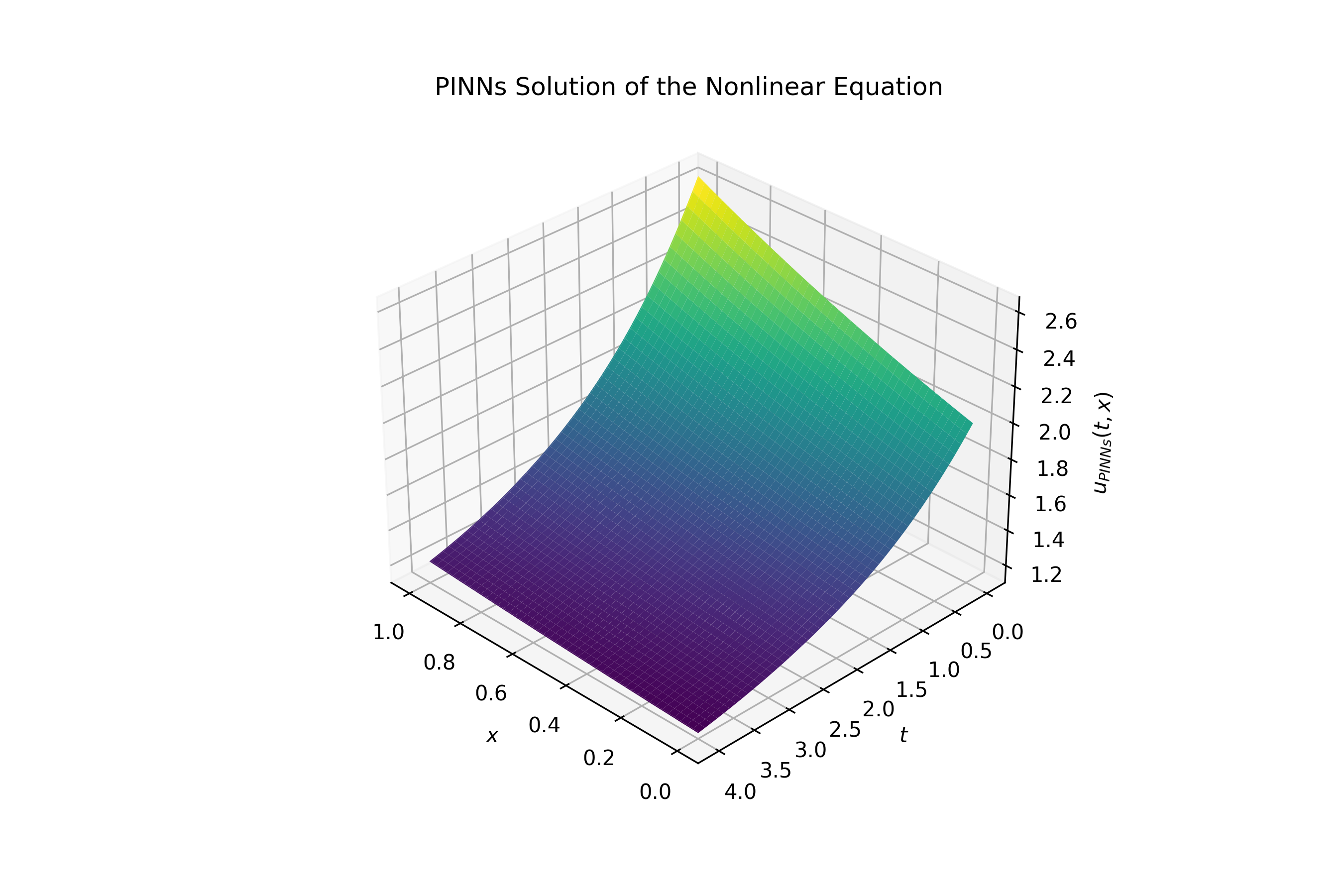}
        \caption{PINNs Solution of the Nonlinear Equation}
        \label{fig:sub3}
    \end{subfigure}
    \begin{subfigure}{0.3\textwidth}
        \centering
        \includegraphics[width=\textwidth]{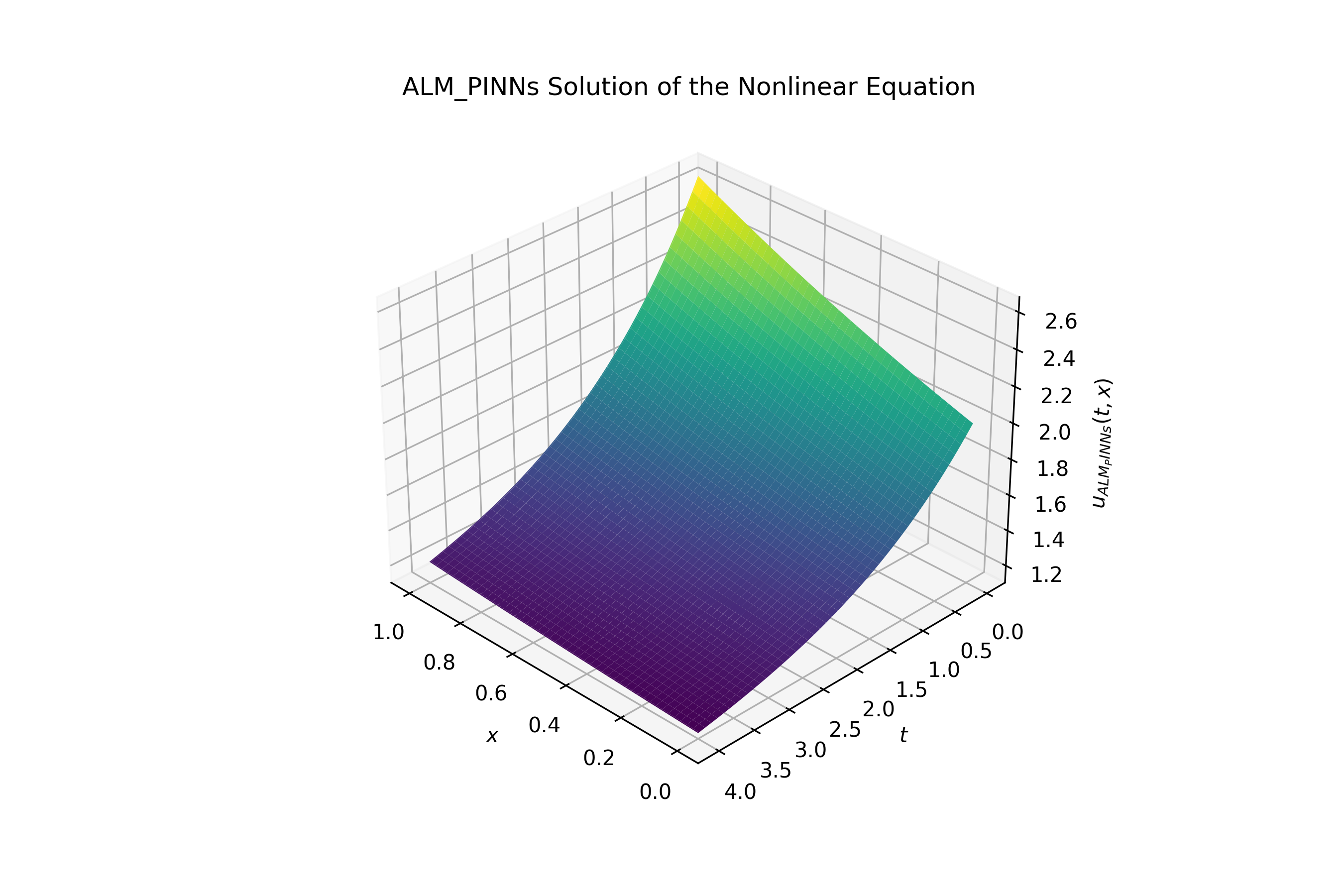}
        \caption{ALM-PINNs Solution of the Nonlinear Equation}
        \label{fig:sub4}
    \end{subfigure}
    
    \begin{subfigure}{0.3\textwidth}
        \centering
        \includegraphics[width=\textwidth]{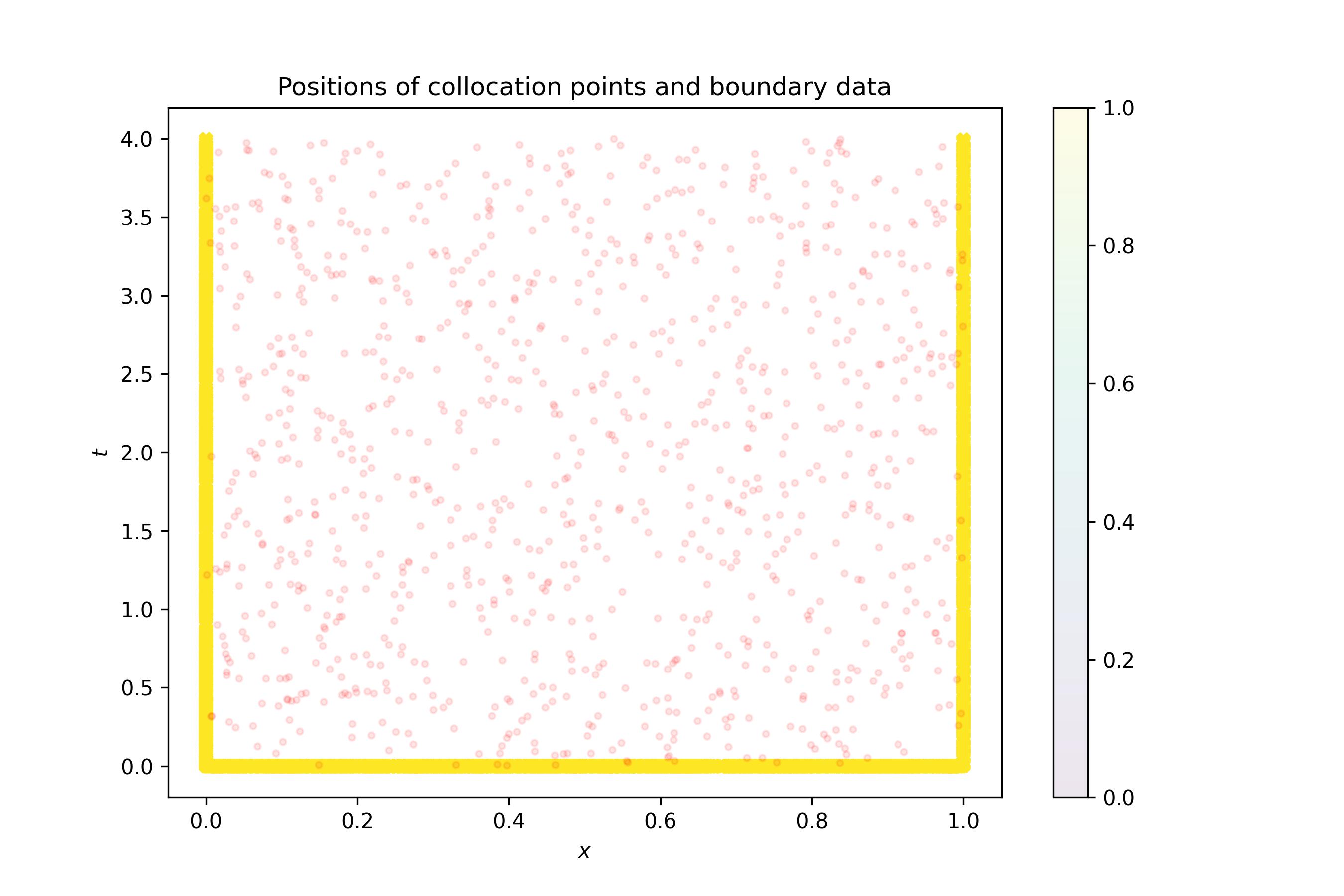}
        \caption{Positions of collocation points and boundary data}
        \label{fig:sub1}
    \end{subfigure}
    \begin{subfigure}{0.3\textwidth}
        \centering
        \includegraphics[width=\textwidth]{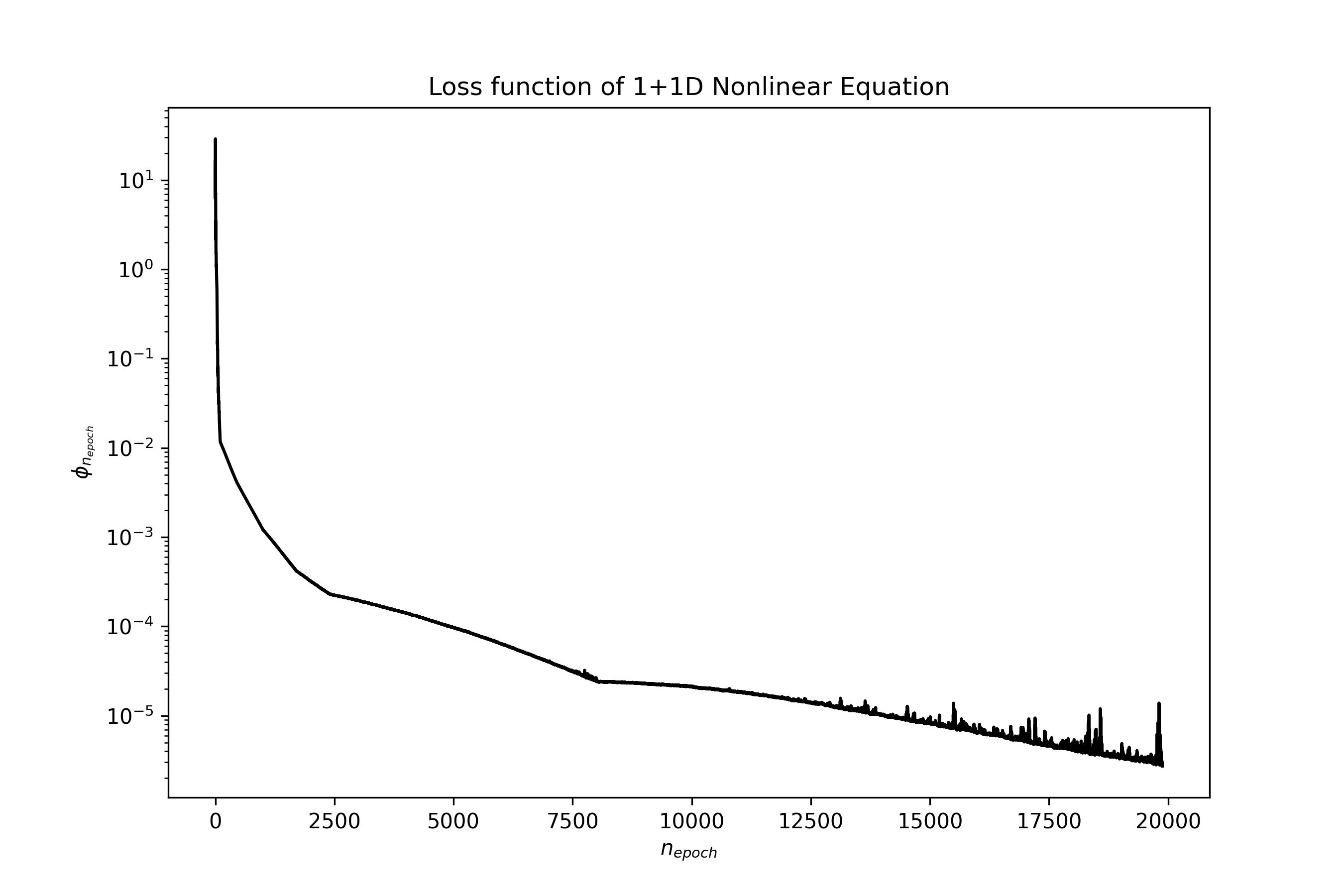}
        \caption{PINNs Evolution of loss function}
        \label{fig:sub5}
    \end{subfigure}
    \begin{subfigure}{0.3\textwidth}
        \centering
        \includegraphics[width=\textwidth]{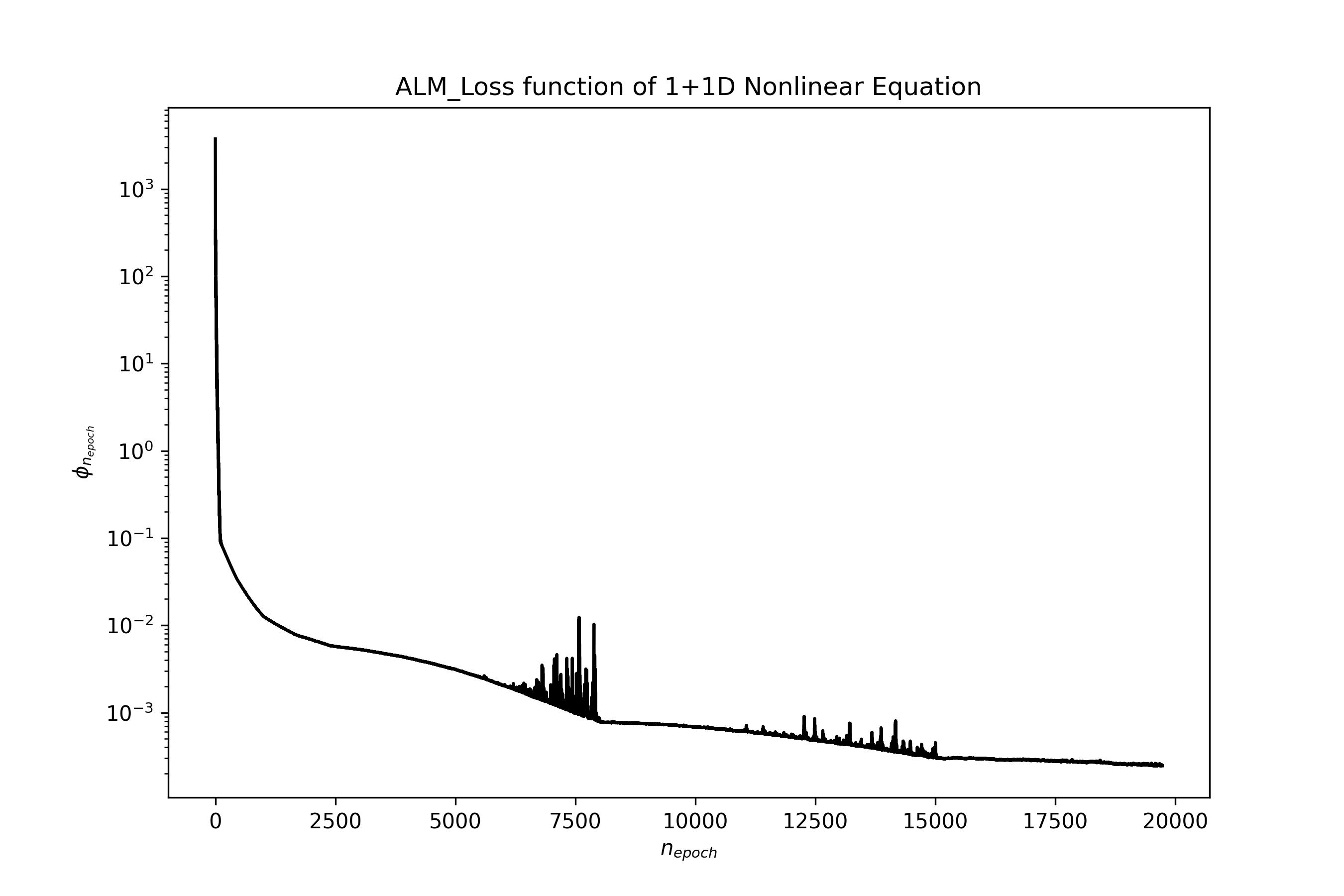}
        \caption{ALM-PINNs Evolution of loss function}
        \label{fig:sub6}
    \end{subfigure}
    
    \caption{Figures of the 1+1 dimensional nonlinear equation when $N_{\mathcal{F}}=1000$,  $N_{\mathcal{B}}=1000$,  $N_{\mathcal{I}}=1000$}
    \label{fig:eg1}
\end{figure}

From the experimental results presented in Table \ref{tab:eg1}, it can be observed that the numerical errors of two algorithms have been computed. The ALM-PINNs achieves more accurate numerical solutions compared to the standard PINNs, thereby validating the main conclusion of Section \ref{sec:ALM-PINNs-forward}. The figures of the sampling points, the loss function, and the solution images are shown in Figure \ref{fig:eg1}.

\begin{table}[htbp]
    \caption{Experimental Parameter Settings in Subsection \ref{eg1}}
    \label{table:compare_1}
    \centering
    \resizebox{\textwidth}{!}{
    \begin{tabular}{cccccc}
        \toprule 
        Hidden Layers & Number of Nodes & Activation Function & Iterations & Optimization Algorithm & Dropout \\
        \midrule
        8 & 40 & tanh & 20000 & Adam & 0.5 \\
        \bottomrule
    \end{tabular}
    }
\end{table}

\begin{table}[h]
    \caption{Results of 1+1 Dimensional Nonlinear Equation Solved by ALM-PINNs and PINNs in Subsection \ref{eg1}}
    \centering
    \begin{tabular}{lcccccc}
        \toprule 
        Models & $\epsilon_r(\hat{u},u)$ & $\epsilon_{\infty}(\hat{u},u)$ & 
        $\epsilon_{a}(\hat{u},u)$ & $N_{\mathcal{F}}$ & $N_{\mathcal{B}}$ & $N_{\mathcal{I}}$\\
        \midrule
         PINNs&1.01693e-05  &3.00999e-03  &5.84584e-07  &1000 &1000 &1000\\
         ALM-PINNs&4.86302e-06  & 1.68999e-03 &1.40627e-07  & 1000 & 1000 &1000 \\
        \bottomrule
    \end{tabular}
    \label{tab:eg1}
\end{table}

\subsection{Viscous Burgers' Equation}
\label{eg2}
Burgers' equation describes the physics of shock wave propagation and reflection. It has important applications in fluid mechanics, nonlinear acoustics and gas dynamics. For a given constant viscosity coefficient $\nu$, the Burgers' equation is also called viscous Burgers' equation. If the viscosity $\nu = 0$, the equation is called inviscid Burgers' equation which governs gas dynamics. 

The nonlinearity of the burgers equation is derived from
$\partial \left(\frac{1}{2}u^2\right)/\partial x$, which is the non-linear term similar in structure to Navier-Stokes equation.

\begin{equation*}
        \frac{\partial u}{\partial t} + \frac{\partial \left(\frac{1}{2}u^2\right)}{\partial x} = \nu \frac{\partial^2 u}{\partial x^2}.
\end{equation*}

Consider the following initial boundary value problem for Burgers' equation
\begin{align*}
    \begin{cases}
        &u_t +  uu_x = \nu u_{xx}, \quad 0 < x < 1, \, t > 0,\\
        &u(x,0) = \sin(\pi x), \quad 0 \leq x \leq 1,\\
        &u(0,t) = 0, \quad u(1,t) = 0, \quad t > 0.
    \end{cases}
\end{align*}

Solving the above problem using Cole-Hopf transformation yields the exact solution of the equation
\begin{equation*}
        u(x,t)=2\pi\nu \frac{\sum_{n=1}^\infty e^{-n^2\pi^2\nu t}nA_n\sin(n\pi x)}{A_0+\sum_{n=1}^\infty e^{-n^2\pi^2\nu t}A_n\cos(n\pi x)},
\end{equation*}
where $A_0,A_n,n=1,2,\cdots$ are Fourier coefficients 
\begin{align*}
        &A_0 = \int_0^1 \exp{(-(2\pi \nu)^{-1}(1-\cos(\pi x)))}\mathrm{d}x,\\
        &A_n = 2\int_0^1 \exp{(-(2\pi \nu)^{-1}(1-\cos(\pi x)))\cos(n\pi x)}\mathrm{d}x.
\end{align*}
\begin{table}[htbp]
        \caption{Experimental Parameter Settings in Subsection \ref{eg2}}
        \label{table:setting2}
        \centering
        \resizebox{\textwidth}{!}{
            \begin{tabular}{cccccc}
                \toprule 
                Hidden Layers & Number of Nodes & Activation Function & Iterations & Optimization Algorithm & Dropout \\
                \midrule
                8 & 40 & tanh & 30000 & Adam & 0.5 \\
                \bottomrule
            \end{tabular}
        }
    \end{table}
    
\begin{table}[htbp]
    \caption{Results of Burgers' equation Solved by ALM-PINNs and PINNs in Subsection \ref{eg2}}
    \centering
    \begin{tabular}{lcccccc}
        \toprule 
        Models & $\epsilon_r(\hat{u},u)$ & $\epsilon_{\infty}(\hat{u},u)$ & 
        $\epsilon_{a}(\hat{u},u)$ & $N_{\mathcal{F}}$ & $N_{\mathcal{B}}$ & $N_{\mathcal{I}}$\\
        \midrule
         PINNs&3.87160e-05 & 1.33683e-02 & 5.02661e-08  &500 &500 &1000\\
         ALM-PINNs&1.35900e-05 & 2.90822e-04 & 6.61664e-09  &500 &500 &1000 \\
        \bottomrule
    \end{tabular}
    \label{tab:eg2}
\end{table}

\begin{figure}[htbp]
    \centering
    \begin{subfigure}{0.3\textwidth}
        \centering
        \includegraphics[width=\textwidth]{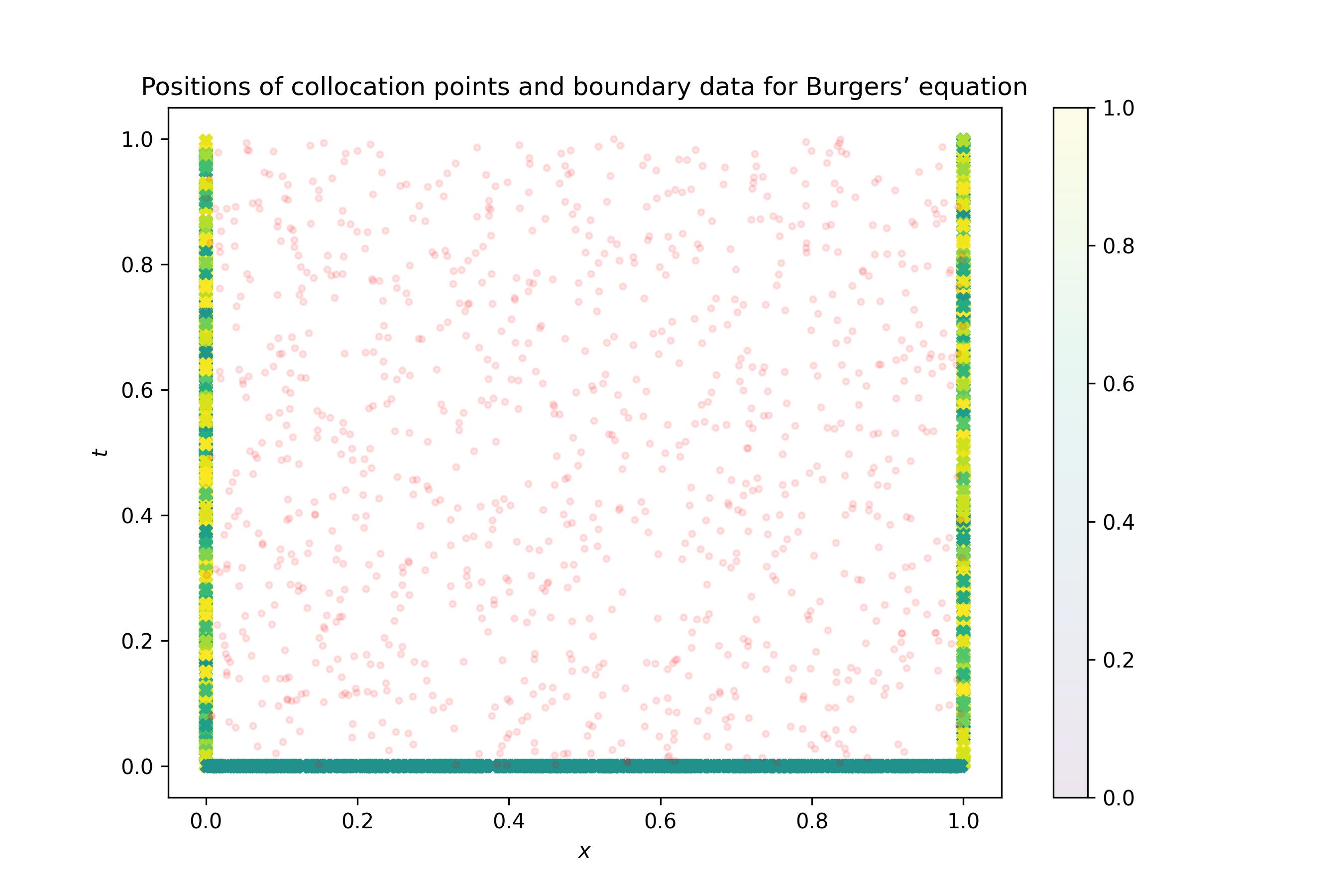}
        \caption{Positions of collocation points and boundary data for Burgers’ equation}
    \end{subfigure}
    \begin{subfigure}{0.3\textwidth}
        \centering
        \includegraphics[width=\textwidth]{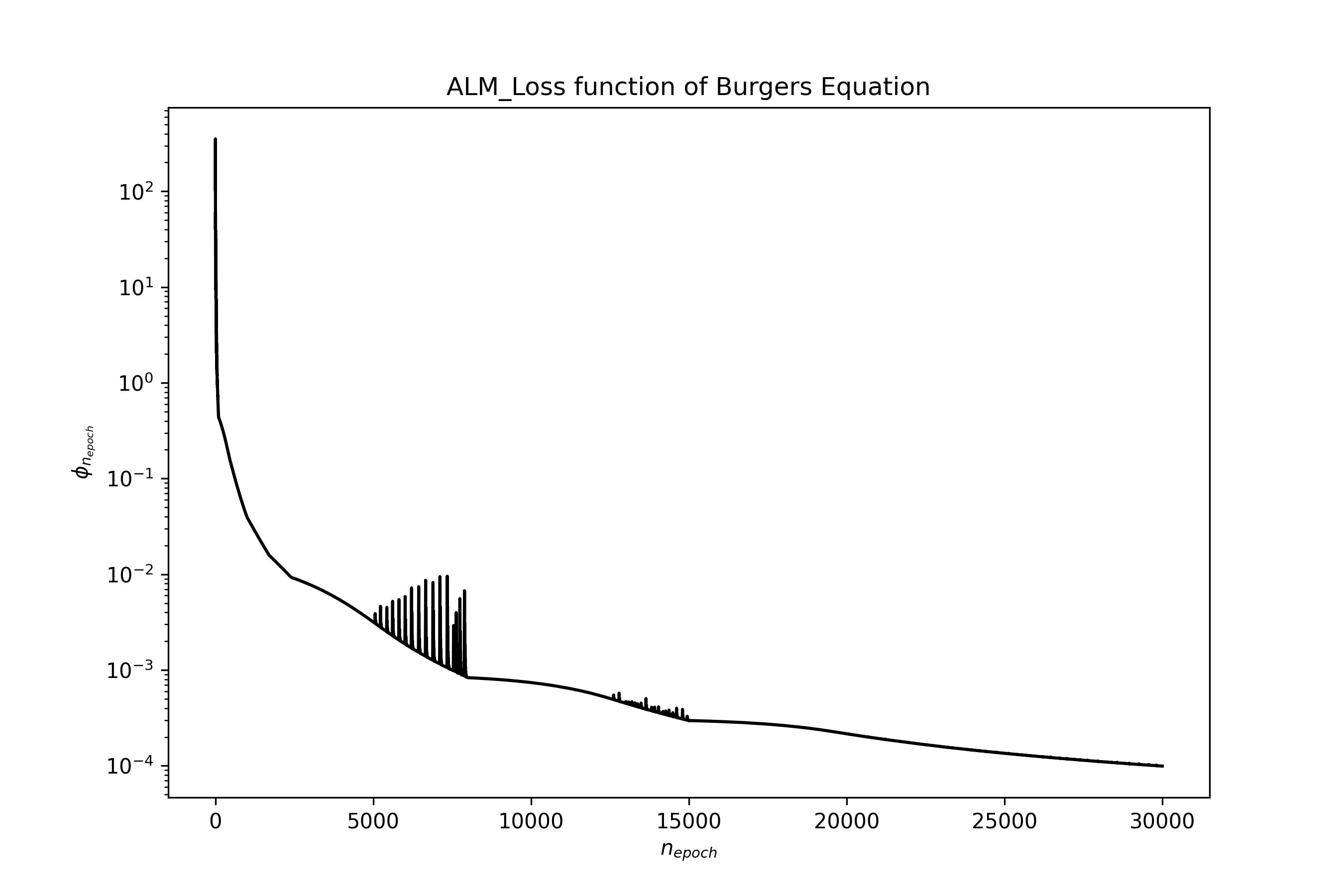}
        \caption{Loss function variation curve of the ALM-PINNs algorithm}
    \end{subfigure}
    \begin{subfigure}{0.3\textwidth}
        \centering
        \includegraphics[width=\textwidth]{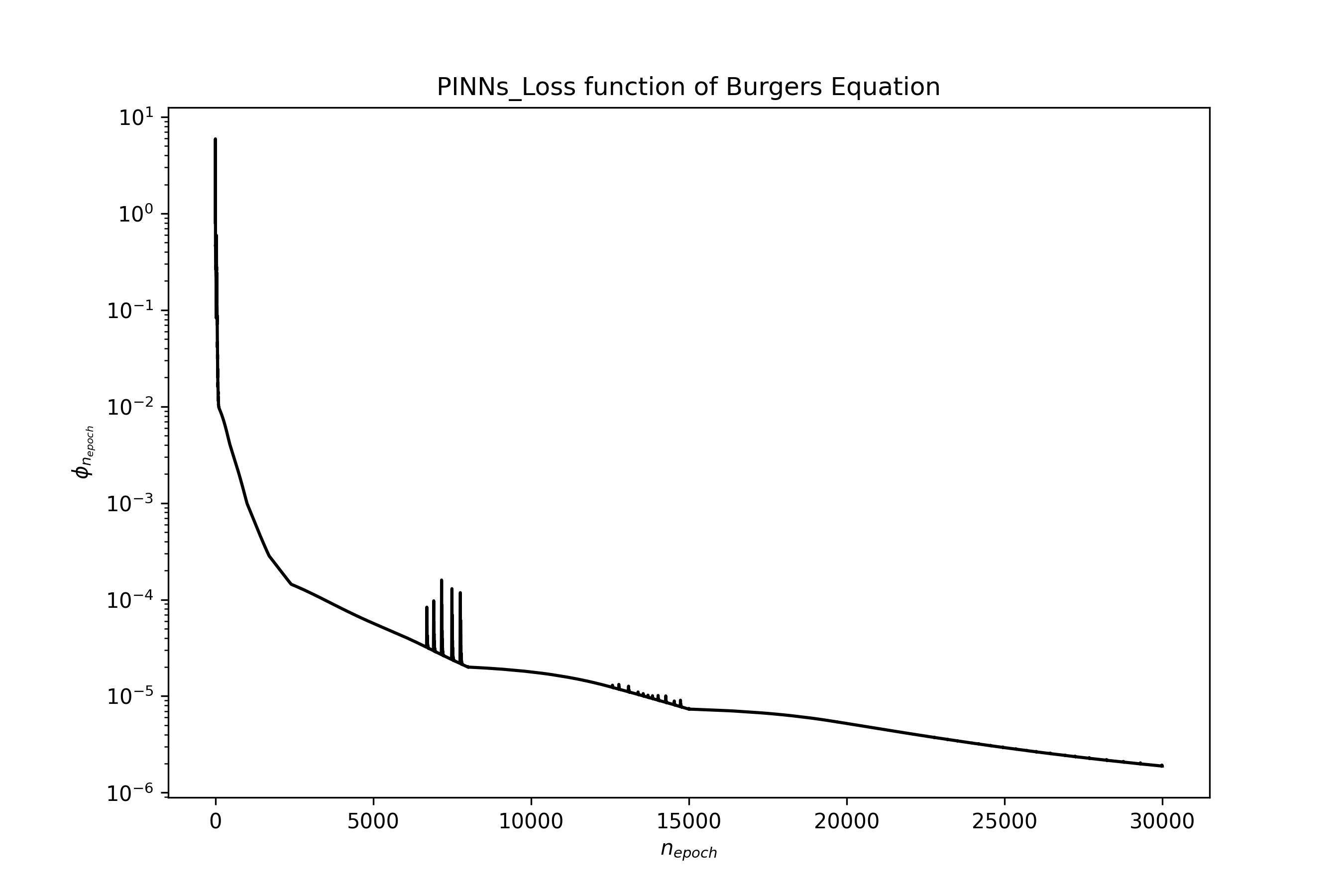}
        \caption{Loss function variation curve of the PINNs algorithm}
    \end{subfigure}
    \begin{subfigure}{0.3\textwidth}
        \centering
        \includegraphics[width=\textwidth]{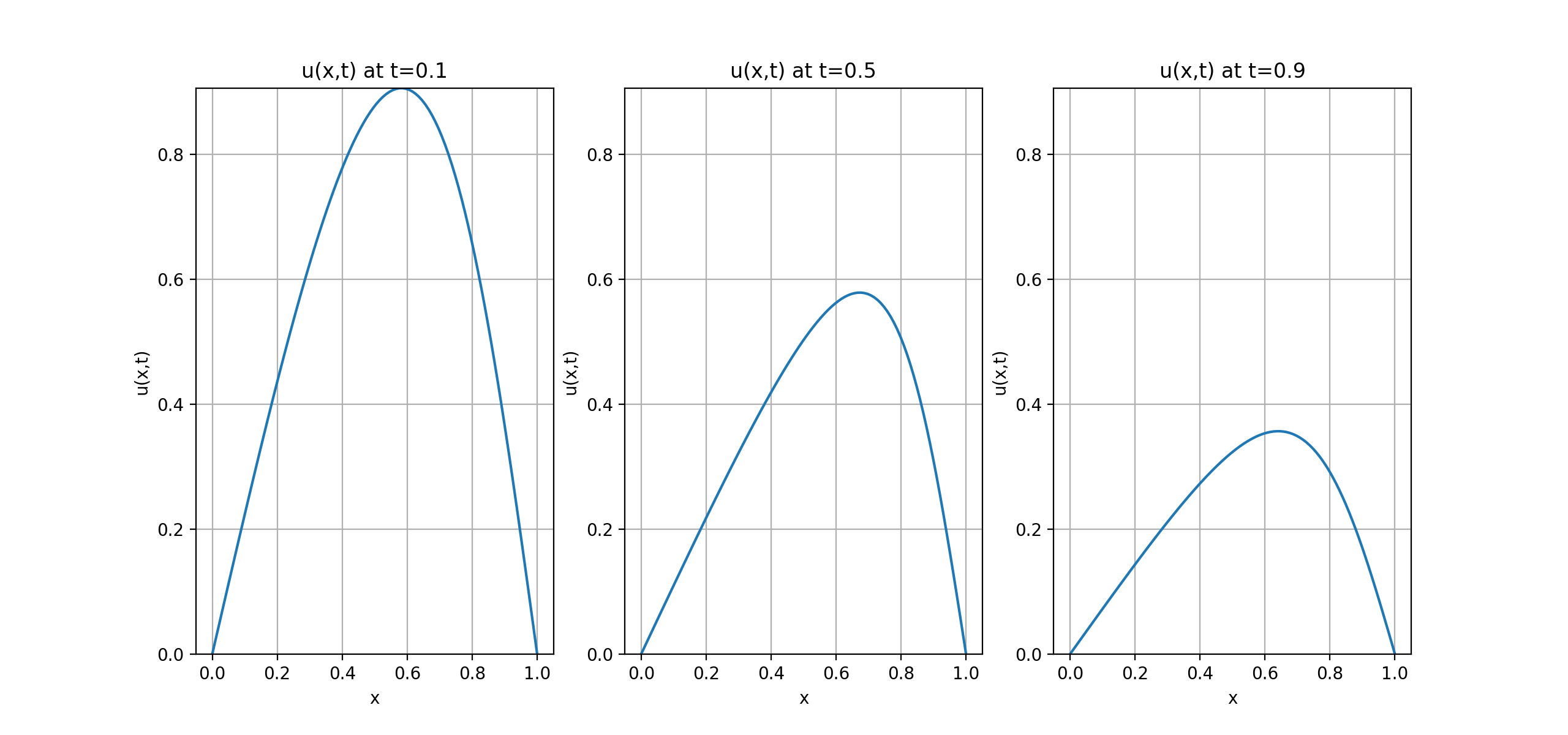}
        \caption{The Burgers' equation plots using the exact solution at $t=0$,$t=0.5$,$t=0.9$}
    \end{subfigure}
    \begin{subfigure}{0.3\textwidth}
        \centering
        \includegraphics[width=\textwidth]{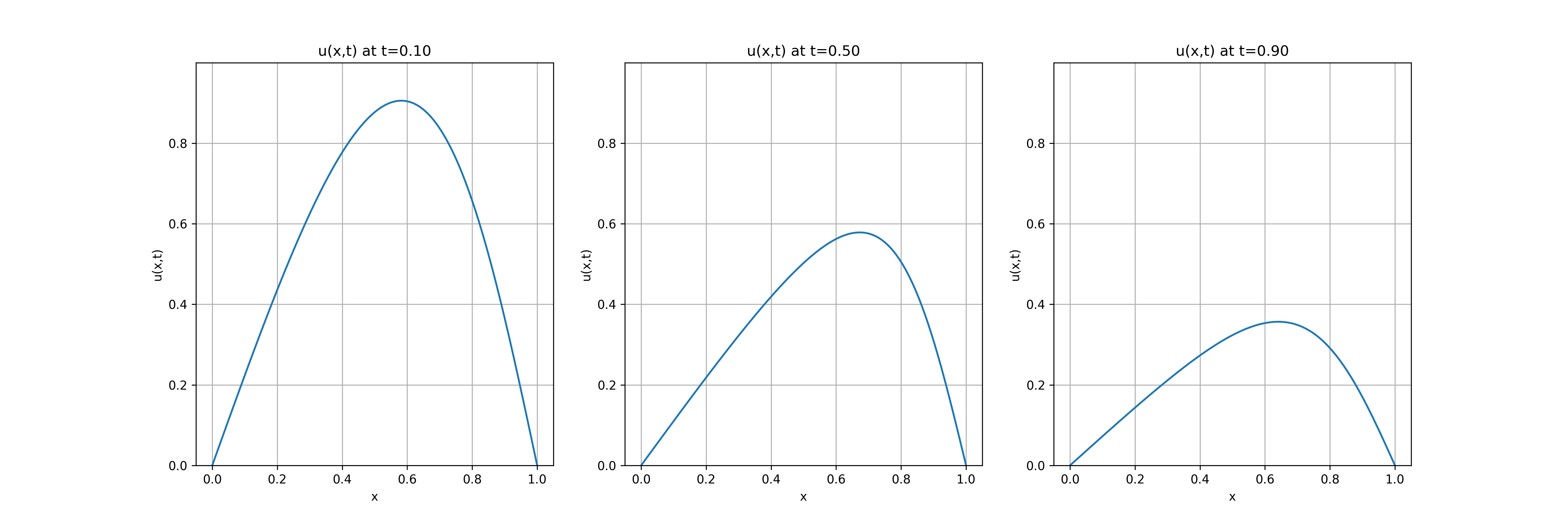}
        \caption{The Burgers' equation plots using ALM-PINNs at $t=0$,$t=0.5$,$t=0.9$}
    \end{subfigure}
    \begin{subfigure}{0.3\textwidth}
        \centering
        \includegraphics[width=\textwidth]{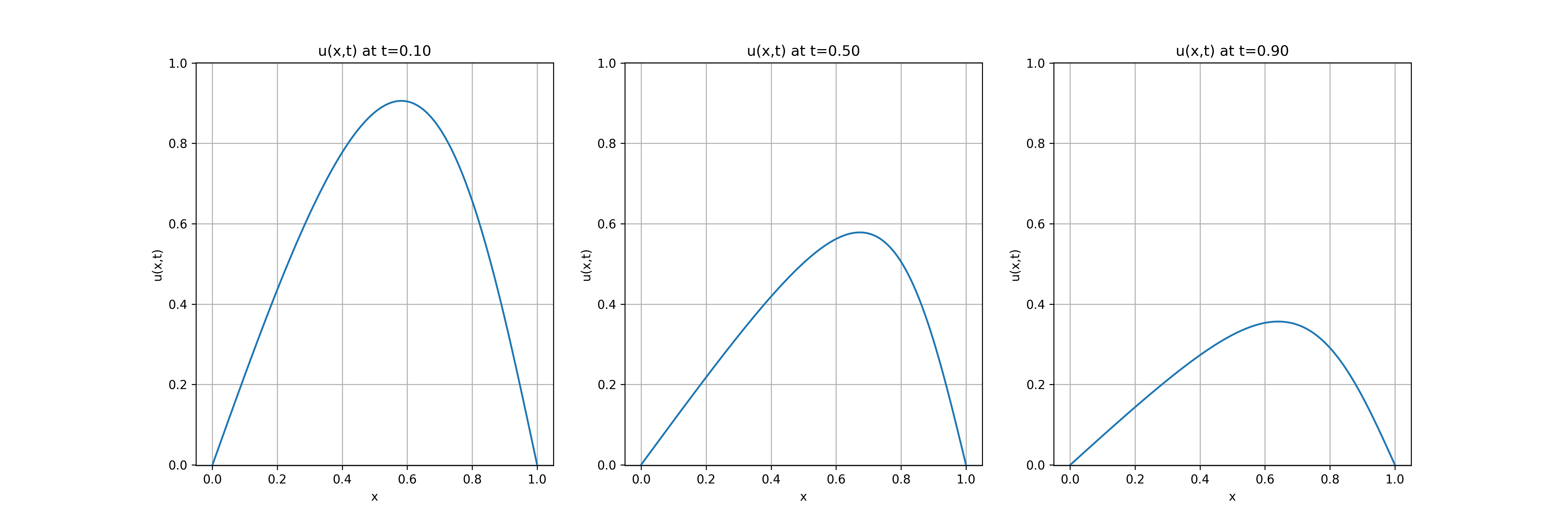}
        \caption{The Burgers' equation plots using PINNs at $t=0$,$t=0.5$,$t=0.9$}
    \end{subfigure}
    \begin{subfigure}{0.3\textwidth}
        \centering
        \includegraphics[width=\textwidth]{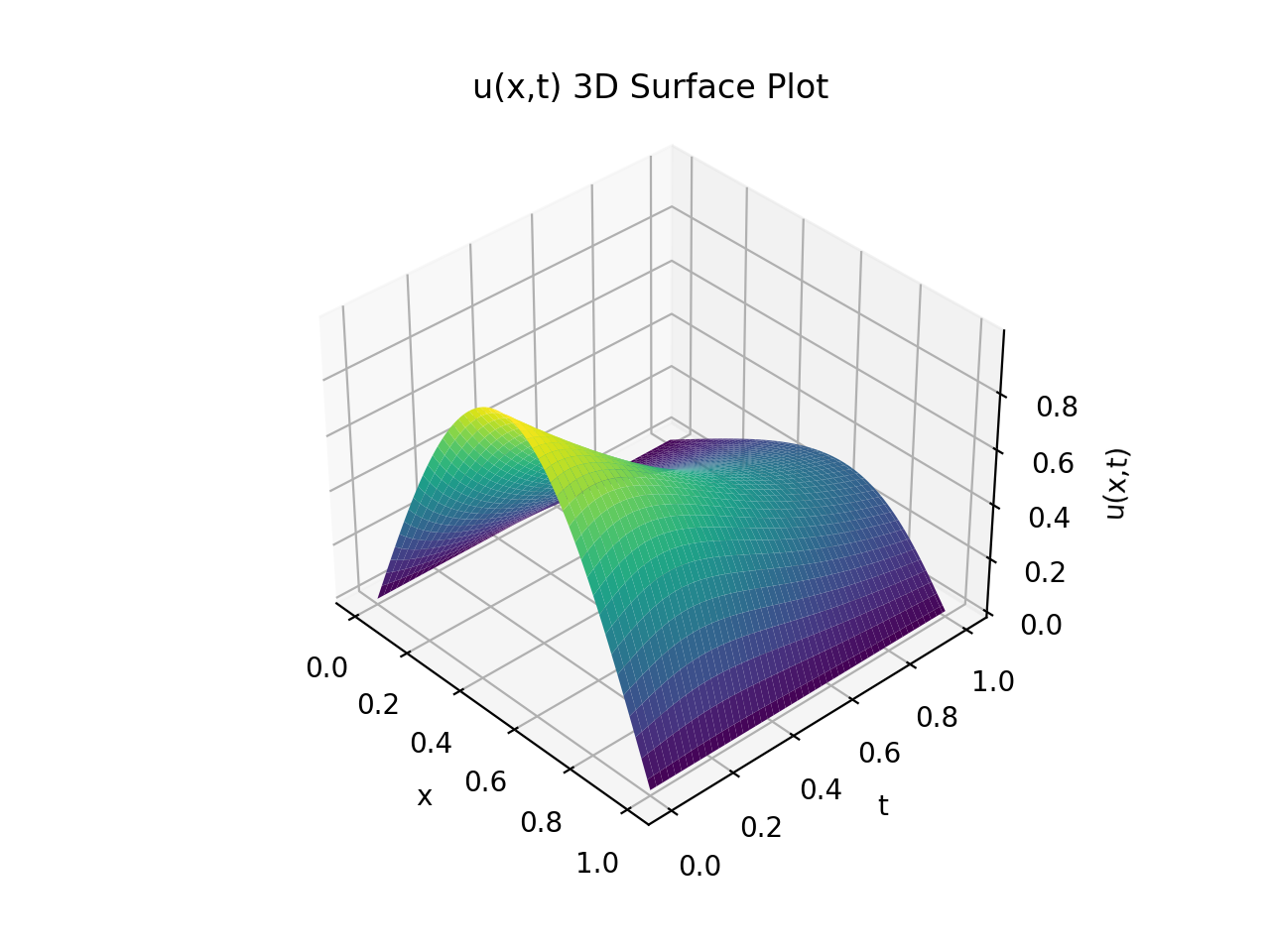}
        \caption{The Burgers' equation plots using the exact solution}
    \end{subfigure}
    \begin{subfigure}{0.3\textwidth}
        \centering
        \includegraphics[width=\textwidth]{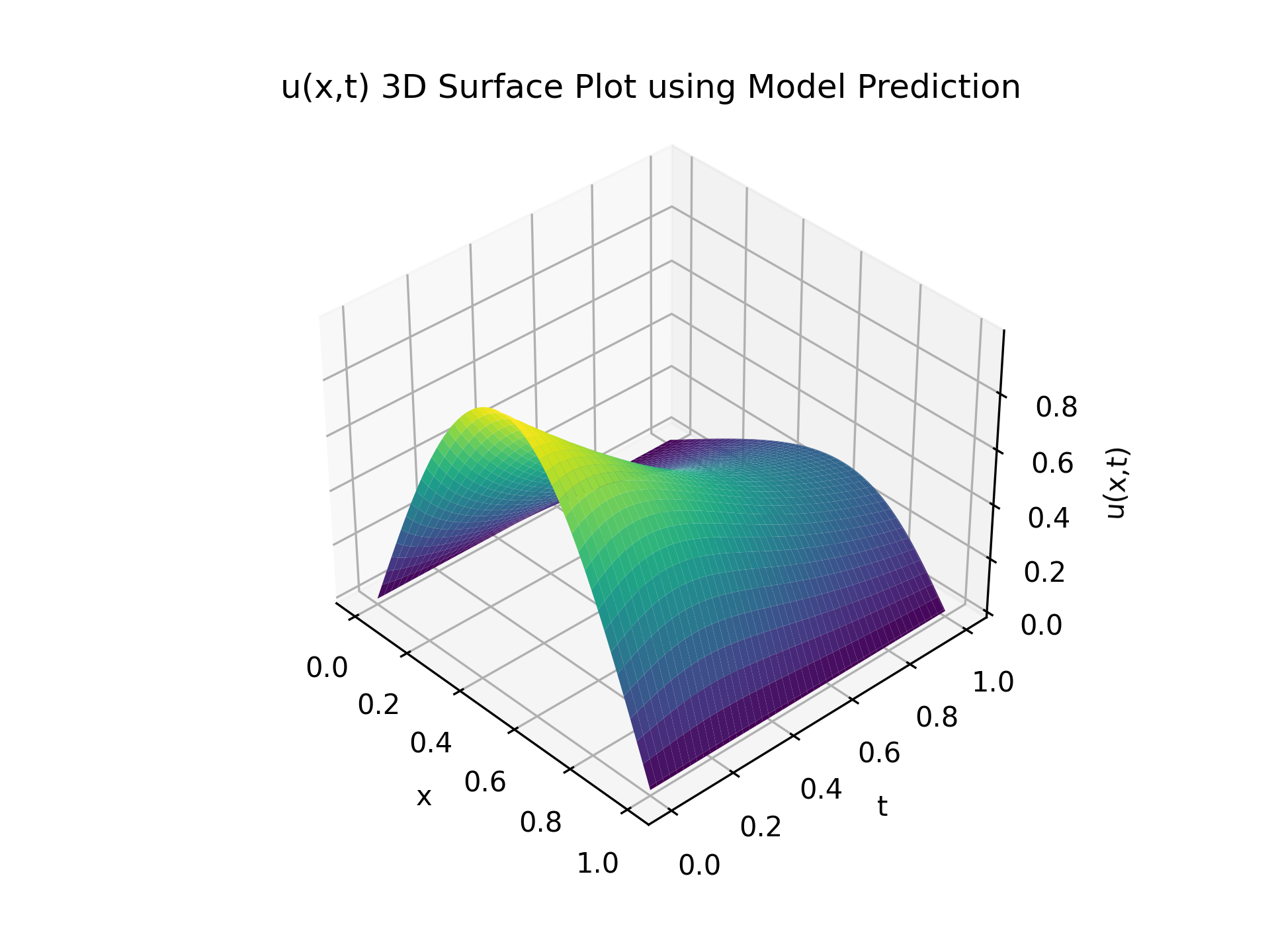}
        \caption{The Burgers' equation plots using ALM-PINNs solution}
    \end{subfigure}
    \begin{subfigure}{0.3\textwidth}
        \centering
        \includegraphics[width=\textwidth]{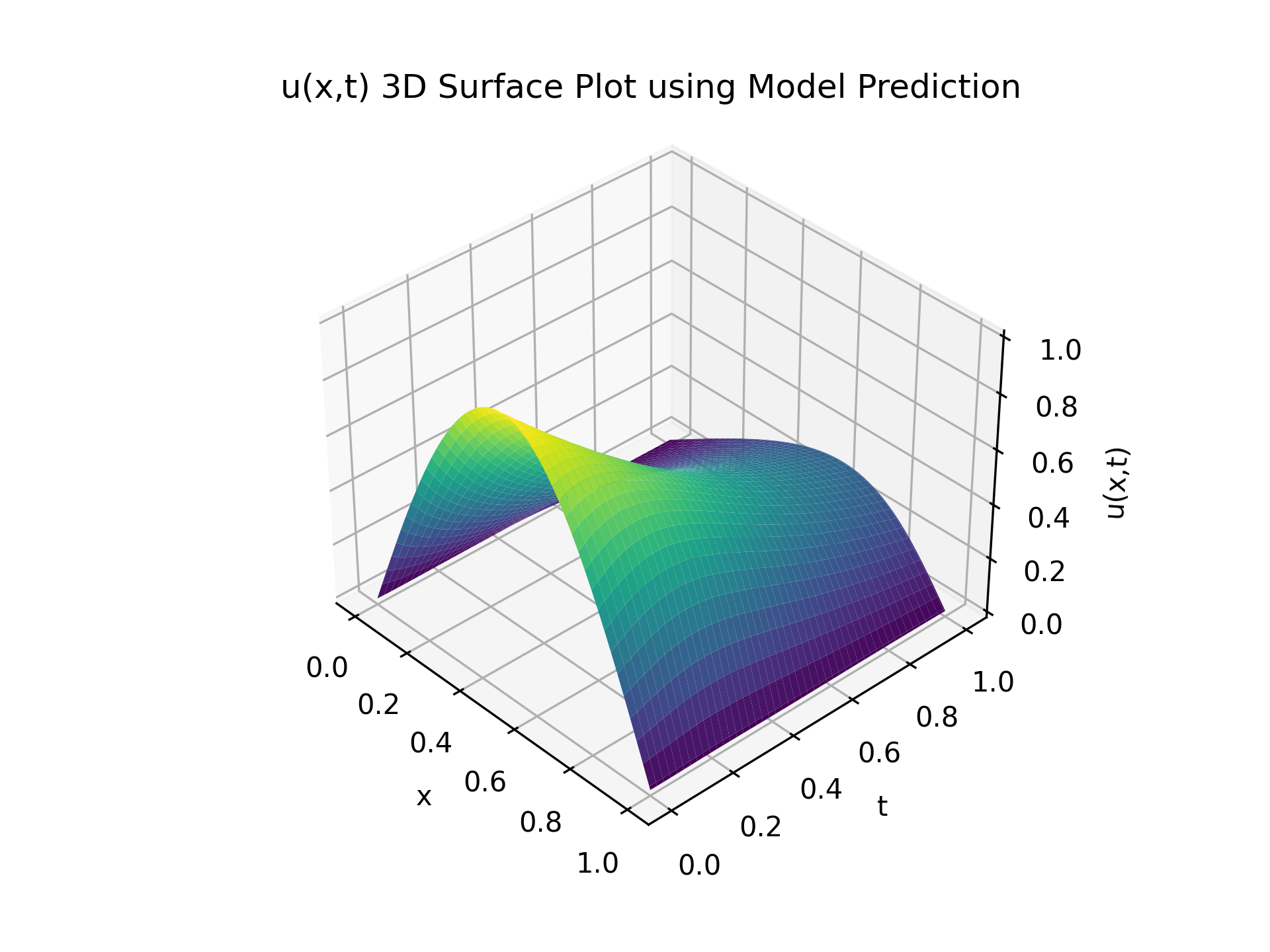}
        \caption{The Burgers' equation plots using PINNs}
    \end{subfigure}
    \caption{Figures of the Burgers' equation when $N_{\mathcal{F}}=500$,  $N_{\mathcal{B}}=500$,  $N_{\mathcal{I}}=1000$}
    \label{inversefig:eg1}
\end{figure}
The experimental parameters for this example are provided in \ref{table:setting2}, and the numerical accuracy of the solutions obtained by the two algorithms is presented in Table \ref{tab:eg2}. The first row of subplots in Figure \ref{inversefig:eg1} shows the sampling points and the loss functions of both algorithms. The second row of subplots presents the slice plots at \(t=0\), \(t=0.5\) and \(t=0.9\), using the exact solution, ALM-PINNs and PINNs algorithms, respectively. The third row of subplots displays the 3D plots generated by the exact solution, ALM-PINNs and PINNs algorithms.

\section{Numerical examples of parameter estimation problems}
\label{sec:inverse}
This section will validate the main conclusions from Section \ref{sec:ALM-PINNs-inverse} through two numerical experiments based on the previous section, as the forward and inverse problems are typically paired.
\subsection{The 1+1 dimensional nonlinear equation}
\label{inverse:eg1}

\begin{equation}
    \left\{
    \begin{array}{lll}
    u_t=\nu_1(u_x)^2+\nu_2uu_{xx}+u-u^2,      &      & t\in [0, 4], x \in [0, 1],\\
    u(x, 0)=\left\{\frac{\tanh(-\frac{x}{4})}{2}+\frac{1}{2}\right\}^{-1},&      & x \in [0, 1],\\
    u(0, t)=\left\{\frac{\tanh(\frac{t}{4})}{2}+\frac{1}{2}\right\}^{-1},   &      & t\in [0, 4],\\
    u(1, t)=\left\{\frac{\tanh(\frac{t}{4}-\frac{1}{4})}{2}+\frac{1}{2}\right\}^{-1},   &      & t\in [0, 4],
    \end{array} \right.  
\end{equation}

\begin{figure}[htbp]
    \centering
    \begin{subfigure}{0.5\textwidth}
        \includegraphics[width=\textwidth]{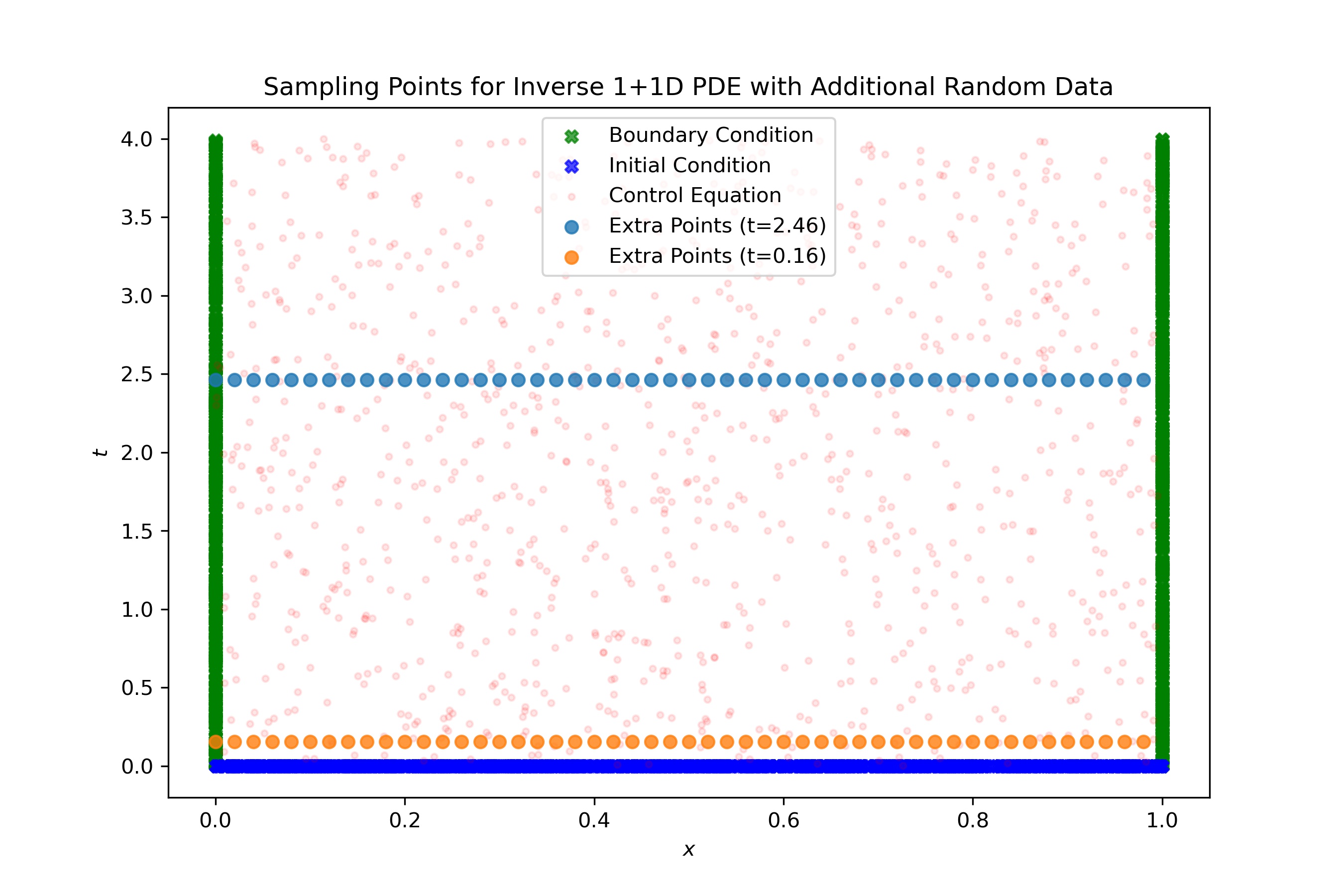}
        \caption{Sampling Points for Inverse 1+1D PDE with Additional Data at $T_i,i=1,2$ Points}
    \end{subfigure}%
    \begin{subfigure}{0.5\textwidth}
        \includegraphics[width=\textwidth]{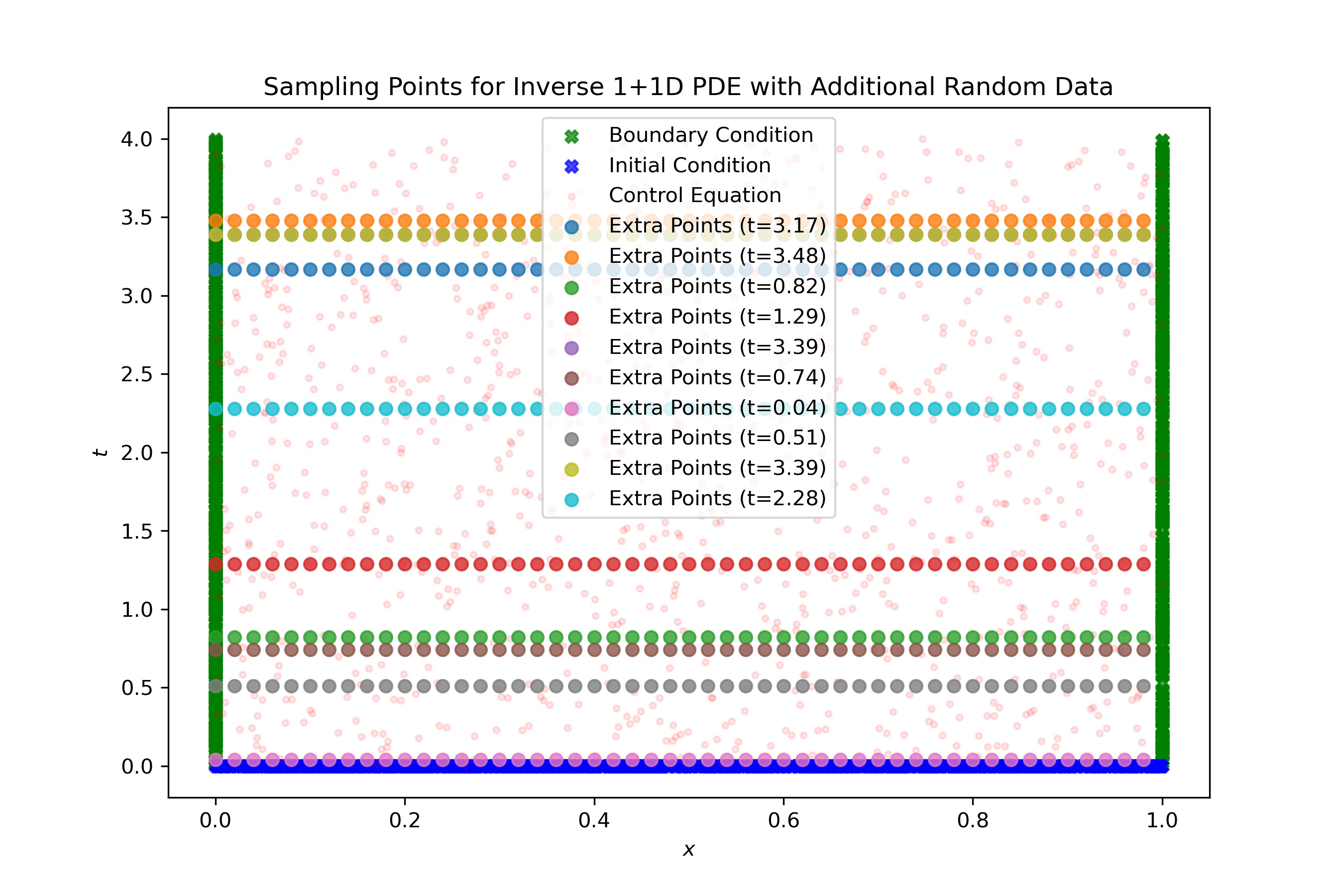}
        \caption{Sampling Points for Inverse 1+1D PDE with Additional Data at $T_i,i=1,\dots,10$ Points}
    \end{subfigure}
    \caption{Schematic of Sampling Points with Different Numbers of Additional Dataset}
    \label{fig:inverse_sample}
\end{figure}

First of all, following the strategy outlined in Section \ref{sec:ALM-PINNs-inverse}, we maintain the same sampling points as those used in the forward problem. Building on this foundation, we introduce varying numbers of additional data points, as illustrated in Figure \ref{fig:inverse_sample}. In Table \ref{inverse_eg1_compare}, we compare the computational results for \(\nu_1\) and \(\nu_2\) obtained using the L2 loss function in the ALM-PINNs algorithm, with an additional dataset of 50 points in \(x\) and varying numbers of points in \(t\) (2, 5, and 10). The final configuration selected for the additional dataset is \(x_{\text{Num}}=50\) and \(t_{\text{Num}}=2\). Table \ref{tab:pinns_alm_combined_1} presents a comparison of the results for \(\nu_1\) and \(\nu_2\) when adding Gaussian errors of 2\% and 20\% while employing L2-PINNs and L2-ALM. Table \ref{tab:pinns_alm_combined_2} compares the results for \(\nu_1\) and \(\nu_2\) when adding Laplace errors of 2\% and 20\% using L2-PINNs, L1-PINNs, L1-ALM, and L2-ALM. An analysis of these results demonstrates that the ALM-PINNs algorithm designed in this study outperforms standard PINNs when utilizing the same error term loss function. This also validates the correctness of Theorem \ref{thm}, highlighting the need to choose different error term loss functions when addressing various types of errors.

\begin{table}[htbp]
\centering
\caption{Comparison of the Relative Errors of \(\nu_1\) and \(\nu_2\) using L2-ALM with 2\% Gaussian noise added}
\label{inverse_eg1_compare}
\begin{tabular}{lccc}
    \toprule
    $x_{Num}=50$ & $t_{Num}=2$ & $t_{Num}=5$ & $t_{Num}=10$ \\
    \midrule
    error\_v\_1 & 8.61585e-03 & 4.93718e-02 & 2.63416e-02 \\
    error\_v\_2 & 1.33055e-03 & 4.76873e-03 & 7.36505e-03 \\
    v\_1 & 2.01723e+00 & 1.90126e+00 & 2.05268e+00 \\
    v\_2 & 1.99734e+00 & 1.99046e+00 & 1.98527e+00 \\
    \bottomrule
\end{tabular}
\end{table}

\begin{table}[htbp]
    \centering
    \caption{Comparison of the results from L2-PINNs and L2-ALM with added Gaussian errors}
    \begin{tabular}{lcc|cc}
        \toprule
        & \multicolumn{2}{c|}{2\% Gaussian error added} & \multicolumn{2}{c}{20\% Gaussian error added} \\
        \midrule
        & L2-PINNs & L2-ALM & L2-PINNs & L2-ALM \\
        \midrule
        error\_v\_1 & 7.9515e-02 & \textbf{8.61585e-03} & 7.6527e-02 & \textbf{9.25946e-03} \\
        error\_v\_2 & 1.9157e-02 & \textbf{1.33055e-03} & 4.7111e-02 & \textbf{3.08764e-03} \\
        v\_1 & 1.8410e+00 & \textbf{2.01723e+00} & 1.8469e+00 & \textbf{2.01852e+00} \\
        v\_2 & 2.0383e+00 & \textbf{1.99734e+00} & 2.0942e+00 & \textbf{1.99382e+00} \\
        \bottomrule
    \end{tabular}
    \label{tab:pinns_alm_combined_1}
\end{table}

\begin{table}[htbp]
    \centering
    \caption{The results in Subsection \ref{inverse:eg1} from L1-PINNs, L2-PINNs, L1-ALM, and L2-ALM include the addition of 2\% and 20\% Laplace relative errors}
    \begin{tabular}{lcccc}
        \toprule
        & L1-PINNs & L2-PINNs & L1-ALM & L2-ALM \\
        \midrule
        \multicolumn{5}{c}{2\% Laplace error added} \\
        \midrule
        error\_v\_1 & 6.2829e-02 & 7.8470e-02 & \textbf{8.8893e-03} & 1.5682e-02 \\
        error\_v\_2 & 1.5658e-02 & 1.6539e-02 & \textbf{2.5761e-03} & 4.1614e-02  \\
        v\_1 & 1.8743e+00 & 1.9737e-02 & \textbf{2.0178e+00} & 1.9686e+00 \\
        v\_2 & 2.0313e+00 & 2.0331e+00 & \textbf{1.9948e+00} & 2.0832e+00  \\
        \midrule
        \multicolumn{5}{c}{20\% Laplace error added} \\
        \midrule
        error\_v\_1 & 4.9946e-01 & 8.8219e-01 & \textbf{1.5276e-02} & 4.5443e-02 \\
        error\_v\_2 & 6.8080e-01 & 2.6266e+00 & \textbf{6.6763e-04} & 2.2302e-03  \\
        v\_1 & 1.4995e+00 & 1.8822e+00 & \textbf{2.0306e+00} & 1.9091e+00 \\
        v\_2 & 1.6808e-01 & 3.6266e-01 & \textbf{1.9987e+00} & 2.0045e+00  \\
        \bottomrule
    \end{tabular}
    \label{tab:pinns_alm_combined_2}
\end{table}

\subsection{Viscous Burgers' Equation}

In this section, the viscosity coefficient $v$ of the inversion equation will be obtained through the exact solution of the example from the perspective of the inverse problem.
\begin{equation*}
    \frac{\partial u}{\partial t} + \nu_1 \frac{\partial \left(\frac{1}{2}u^2\right)}{\partial x} = \nu_2 \frac{\partial^2 u}{\partial x^2}.
\end{equation*}
\begin{align*}
\begin{cases}
    &u_t +  \nu_1 uu_x = \nu_2 u_{xx}, \quad 0 < x < 1, \, t > 0,\\
    &u(x,0) = \sin(\pi x), \quad 0 \leq x \leq 1,\\
    &u(0,t) = 0, \quad u(1,t) = 0, \quad t > 0.
\end{cases}
\end{align*}
We conducted a total of 20 small experiments across 3 major groups, each classified based on different noise distributions, to validate the effectiveness of the proposed Algorithm \ref{thm}.

\textbf{First group of experiments}: The noise was assumed to follow a normal distribution with relative noise levels of 2$\%$ and 20$\%$. Without and with the regularization strategies proposed in Section \label{ill-posed}, we performed parameter inversion experiments for the Burgers' equation using the ALM-PINNs-L2 loss function and the PINNs-L2 loss function, respectively. The results showed that under the regularization strategy, the ALM-PINNs-L2 loss function produced the optimal parameter inversion results.

\textbf{Second group of experiments}: The noise was assumed to follow a Laplace distribution. Similarly, we conducted experiments at relative noise levels of 2$\%$ and 20$\%$, using the ALM-PINNs-L2 loss function, ALM-PINNs-L1 loss function, PINNs-L2 loss function, and PINNs-L1 loss function, respectively. We also repeated the experiments with the regularization strategy proposed in Section 4. The results showed that under the regularization strategy, the ALM-PINNs-L1 loss function provided the most accurate parameter inversion results.

\textbf{Third group of experiments}: The noise was assumed to follow a log-normal distribution. As in the previous experiments, we tested at relative noise levels of 2$\%$ and 20$\%$, using the ALM-PINNs-L2 loss function, ALM-PINNs-LogN loss function, PINNs-L2 loss function, and PINNs-LogN loss function. Regularization strategies were applied again for repeated experiments. The results demonstrated that, under the regularization strategies proposed in Section \label{ill-posed}, the ALM-PINNs-LogN loss function produced the most optimal parameter inversion results.
\begin{table}[htbp]
    \centering
    \caption{Results from L2-PINNs and L2-ALM with Gauss relative errors of 2\% and 20\%}
    \begin{tabular}{lcc|cc}
        \toprule
        & \multicolumn{2}{c|}{2\% Gaussian error added} & \multicolumn{2}{c}{20\% Gaussian error added} \\
        \midrule
        & L2-PINNs & L2-ALM & L2-PINNs & L2-ALM \\
        \midrule
        error\_v\_1 & 4.3660e-01 & \textbf{6.3586e-04} & 6.5463e-01 & \textbf{2.7018e-02} \\
        error\_v\_2 & 4.1259e-01 & \textbf{9.9999e-02} & 3.7794e-01 & \textbf{1.6421e-01} \\
        v\_1 & 1.4366e+00 & \textbf{1.0006e+00} & 1.6546e+00 & \textbf{1.0270e+00} \\
        v\_2 & 1.4126e-01 & \textbf{9.0000e-02} & 1.3779e-01 & \textbf{8.3578e-02} \\
        \bottomrule
    \end{tabular}
\end{table}


\begin{table}[htbp]
    \centering
    \caption{Results from L1-PINNs, L2-PINNs, L1-ALM, and L2-ALM. The table presents results with a Laplace relative error of 2\% and 20\%.}
    \begin{tabular}{lcccc}
        \toprule
        & L1-PINNs & L2-PINNs & L1-ALM & L2-ALM \\
        \midrule
        \multicolumn{5}{c}{2\% Laplace error added} \\
        \midrule
        error\_v\_1 & 4.9946e-01 & 8.8219e-01 & \textbf{6.1028e-03} & 2.7731e-02 \\
        error\_v\_2 & 6.8080e-01 & 2.6266e+00 & \textbf{7.9999e-02} & 1.5904e-01  \\
        v\_1 & 1.4995e+00 & 1.8822e+00 & \textbf{1.0061e+00} & 1.0277e+00 \\
        v\_2 & 1.6808e-01 & 3.6266e-01 & \textbf{9.2000e-02} & 8.4096e-02 \\
        \midrule
        \multicolumn{5}{c}{20\% Laplace error added} \\
        \midrule
        error\_v\_1 & 9.3719e-01 & 1.0145e+00 & \textbf{2.7149e-02} & \textbf{3.3715e-02} \\
        error\_v\_2 & 3.3421e+00 & 6.9220e+00 & \textbf{1.0246e-01} & \textbf{2.0000e-01}  \\
        v\_1 & 1.9372e+00 & 2.0145e+00 & \textbf{1.0271e+00} & \textbf{9.6628e-01} \\
        v\_2 & 4.3421e-01 & 7.9220e-01 & \textbf{8.9754e-02} & \textbf{7.9999e-02}  \\
        \bottomrule
    \end{tabular}
    \label{tab:pinns_alm_combined}
\end{table}

\begin{table}[htbp]
    \centering
    \caption{Results from Log-PINNs, L2-PINNs, Log-ALM, and L2-ALM with Log-normal errors of 2\% and 20\% added.}
    \begin{tabular}{lcccc}
        \toprule
         & Log-PINNs & L2-PINNs & \textbf{Log-ALM} & L2-ALM \\
        \midrule
        \multicolumn{5}{c}{2\% Log-Normal error added} \\
        \midrule
         error\_v\_1 & 2.1556e-02 & 1.0676e-01 & \textbf{5.5206e-03} & 1.8411e-02 \\
        error\_v\_2 & 1.8000e-01 & 2.0000e-01 & \textbf{6.6001e-02} & 8.5856e-02 \\
        v\_1 & 9.7844e-01 & 8.9324e-01 & \textbf{1.0055e+00} & 1.0184e+00 \\
        v\_2 & 1.1800e-01 & 1.2000e-01 & \textbf{9.3399e-02} & 9.1414e-02 \\
        \midrule
        \multicolumn{5}{c}{20\% Log-Normal error added} \\
        \midrule
        error\_v\_1 & 5.6458e-01 & 8.7570e-01 & \textbf{3.2143e-02} & 2.3169e-02 \\
        error\_v\_2 & 1.8749e-01 & 1.9362e-01 & \textbf{3.6427e-02} & 8.4000e-02 \\
        v\_1 & 4.3542e-01 & 1.2430e-01 & \textbf{9.6786e-01} & 9.7683e-01 \\
        v\_2 & 1.1875e-01 & 1.1936e-01 & \textbf{9.6357e-02} & 9.1600e-02 \\
        \bottomrule
    \end{tabular}
\end{table}

\section{Conclusion}
\label{sec:conclusion}
The primary objective of this paper is to design the ALM-PINNs algorithm for solving nonlinear PDEs and addressing parameter inversion problems. The key research contributions are as follows:

In the context of solving PDEs, the ALM-PINNs algorithm was developed and successfully applied to two differential equation cases, yielding highly accurate solutions. For parameter inversion, a distinct optimization problem, different from the forward problem, was proposed. A construction theorem for the loss function on additional datasets was derived, and strategies were introduced to overcome the ill-posedness in neural networks. Based on the results of the forward problem, parameter inversion was carried out, confirming the effectiveness of the proposed algorithm.

The experimental results from the selected numerical examples demonstrate that the ALM-PINNs algorithm not only solves various types of PDEs but also achieves higher numerical accuracy compared to the standard PINNs.

With the continuous innovation in computer technology, the computational capabilities of artificial intelligence algorithms are gradually gaining recognition in the field of mathematical physics. PINNs have been a research hotspot in recent years, and while this paper provides a preliminary investigation into the PINNs algorithm, the following valuable research directions are suggested for future exploration:

(1). The field of Physical Informed Machine Learning (PIML), which includes branches such as PINNs, is a vast area. Incorporating physical prior information into algorithms is a research hotspot in this field. Karniadakis et al. have proposed several future research directions. Future research could explore the widespread application of artificial intelligence in applied mathematics.

(2). The theoretical foundation and convergence analysis of PINNs and its derivative algorithms remain areas for further research. Additionally, the criteria for selecting the size of supplementary datasets also require further investigation.

(3). Complex equations often lack analytical solutions, thus requiring higher precision in numerical solutions. Therefore, artificial intelligence algorithms still need theoretical breakthroughs to meet the demands of high-precision numerical solutions.

\section*{Acknowledgements}
This research is partially supported by the National Natural Science Foundation of China (Grant No. 12371428 and 11871435) and the Fundamental Research Funds for the Central Universities (Grant No. CXJJ-2024-443).


\bibliographystyle{elsarticle-num}
\bibliography{elsarticle-template-num}  

\end{document}